\documentclass[11pt]{amsart}
\usepackage{float}
\usepackage{textgreek}
\usepackage[a4paper,margin=1in]{geometry}
\usepackage{amsmath,amssymb,amsthm,mathtools}
\usepackage{enumitem}
\usepackage{lipsum}
\usepackage{microtype}
\usepackage{hyperref}
\usepackage{amsthm}
\usepackage{soul}
\usepackage{enumitem}
\usepackage[numbers]{natbib}
\usepackage{mathtools}
\usepackage{tikz-cd}
\usepackage{bigints}
\usepackage{xfrac}
\usepackage{mathrsfs}
\usepackage{comment}
\usepackage{verbatim}
\usepackage{adjustbox}
\usepackage{graphicx}
\usepackage{import}
\usepackage{upgreek}  
\hypersetup{colorlinks=true,linkcolor=blue,citecolor=blue,urlcolor=blue}
\newtheorem{theorem}{Theorem}[section]
\newtheorem{proposition}[theorem]{Proposition}
\newtheorem{lemma}[theorem]{Lemma}
\newtheorem{corollary}[theorem]{Corollary}
\newtheorem{definition}[theorem]{Definition}
\newtheorem{remark}[theorem]{Remark}
\newtheorem*{lemma*}{Lemma}

\newcommand{\R}{\mathbb{R}}

\newcommand{\N}{\mathbb{N}}

\newcommand{\eps}{\varepsilon}
\newcommand{\ang}[1]{\langle #1\rangle}
\newcommand{\op}{\operatorname{Op}}

\newcommand{\re}{\operatorname{Re}}

\newcommand{\calU}{\mathcal U}

\newcommand{\Lip}{\operatorname{Lip}}
\newcommand{\vertiii}[1]{{\left\vert\kern-0.25ex\left\vert\kern-0.25ex\left\vert #1 
    \right\vert\kern-0.25ex\right\vert\kern-0.25ex\right\vert}}
    \usepackage{xcolor}
    \usepackage{hyperref}
         \hypersetup{colorlinks, citecolor=blue, filecolor=blue, linkcolor=black, urlcolor=blue}
    \usepackage{soul}
\usepackage{enumitem}
    \usepackage{amsmath,amssymb,amsthm,mathtools}
    \usepackage{tikz-cd}
    \usepackage{bigints}
    \usepackage{xfrac}
    \usepackage{mathrsfs}
    \usepackage{comment}
    \usepackage{verbatim}
    \usepackage{adjustbox}
    \usepackage{graphicx}
    \usepackage{import}
     \usepackage{upgreek}

\newtheorem{mainthm}{Theorem}

\newcommand{\Leb}{\operatorname{Leb}}
\newcommand{\loc}{\mathrm{loc}}
\newcommand{\intt}{\operatorname{int}}

\usepackage{marginnote}

\numberwithin{equation}{section}

\newtheorem*{proposition*}{Proposition}

\newtheorem*{question*}{Question}
\newtheorem*{theorem*}{Theorem}
\newtheorem*{claim*}{Claim}

\theoremstyle{definition}

\newtheorem{example}[theorem]{Example}

\theoremstyle{remark}

\makeatletter
\renewcommand\NAT@bibsetnum[1]{%
  \settowidth\labelwidth{\@biblabel{#1}}%
  \setlength{\leftmargin}{2mm}%
  \setlength{\rightmargin}{0pt}%
  \setlength{\itemindent}{0pt}%
  \setlength{\listparindent}{0pt}%
  \setlength{\parsep}{0pt}%
  \setlength{\itemsep}{\bibsep}%
  \setlength{\labelsep}{0.5em}%
}
\makeatother
\begin{document}
\title[Generalized Horseshoe Maps]{Generalized Horseshoe Maps part I:\\limit laws and absolute continuity of physical measure
}
\author{Abbas Fakhari$\,^1$, Maisam Hedyehloo$\,^2$ and Mohammad Soufi$\,^3$} 
\address{$^1\,$Department of Mathematics,
Shahid Beheshti University, 19839 Tehran, Iran.\it{a\_fakhari@sbu.ac.ir}}
\address{$^2\,$School of Mathematics, Institute for Research in Fundamental Sciences (IPM), 19395-5746, Tehran, Iran.
\it{maisam.hedyehloo@ipm.ir}}
\address{$^3\,$Instituto de Matemática e Estatística, Universidade do Estado do Rio de Janeiro, Rio de Janeiro, Brazil.
\it{mohammad@ime.uerj.br}}
\date{} 
\begin{abstract}
We apply thermodynamic formalism to a generalized horseshoe map. We prove that a tailored anisotropic Banach space 
with weighted norms yields a spectral gap for the transfer operator, implying the existence of a unique physical measure and 
limit theorems. Under the virtually expanding condition, this measure is absolutely continuous with respect to Lebesgue measure, with density in the Sobolev space
$H^\mu$, for some $\mu<1/2$.
\end{abstract}

\maketitle
\tableofcontents
\section{Introduction}

A central aim of ergodic theory is to understand the long-term statistical
behavior of typical orbits of a dynamical system. Rather than following individual trajectories, one studies the asymptotic behaviour of empirical measures and their convergence to invariant probability measures, which describe the statistical distribution of almost every initial condition. 
For applications, however, one would like this statistical description to
hold for a set of initial conditions that is observable with respect to
Lebesgue measure.  This leads to the notion of a physical measure and, in
hyperbolic dynamics, to Sinai--Ruelle--Bowen measures.

For smooth hyperbolic systems, SRB measures are characterized by absolute
continuity of their conditional measures along unstable manifolds and,
under the standard hypotheses, describe the asymptotic statistics of a
positive-Lebesgue-measure set of initial conditions.  Ruelle's construction
for Axiom~A attractors is one of the foundational results in this direction
\cite{R76}; we refer to Young's survey \cite{Y02} for an account of
the physical interpretation of SRB measures and of the classes of systems
for which they are known to exist.  These ideas extend beyond globally
smooth dynamics to systems that are hyperbolic only branchwise and possess
discontinuities, singularities, or countably many smoothness domains.

Generalized horseshoe maps form an important class of such systems.
Jakobson and Newhouse introduced them as a two-dimensional counterpart of
the one-dimensional folklore theorem for piecewise expanding maps
\cite{JN96} and as a typical model for the existence of physical measures.  
Their model consists of a countable collection of full-height
vertical strips in a planar unite square on which the map is piecewise \(C^2\), 
uniformly hyperbolic, and
allowed to have unbounded first derivatives.  Under suitable geometric,
hyperbolicity, and distortion assumptions, Jakobson--Newhouse developed
the corresponding stable and unstable manifold theory and constructed an
ergodic SRB measure \cite{JN00}.  S\'anchez-Salas subsequently developed an
SRB-measure construction in a broader piecewise-hyperbolic setting
\cite{SS01}.  More recently in \cite{FaKhSo24}, it was proved the
existence of an invariant measure absolutely continuous with respect to
planar Lebesgue measure for generalized horseshoe maps satisfying
additional fatness and transversality assumptions.

The existence of an SRB or physical measure provides a first statistical
description of the dynamics, but it naturally leads to finer questions.
One would like to determine uniqueness and ergodicity, rates of convergence
to equilibrium, decay of correlations, spectral properties of the
associated transfer operator, and probabilistic limit laws for Birkhoff
sums.  When an invariant measure is absolutely continuous with respect to
ambient Lebesgue measure, the regularity of its density is another natural
problem.  In dissipative hyperbolic systems, however, the invariant measure
is generally singular in the ambient space; the appropriate replacement is
to study its regularity as an element of a distribution space adapted to
the stable and unstable directions.

Finer statistical properties of generalized horseshoe maps were first
investigated by Jakobson through countable symbolic models.  Using
thermodynamic formalism for the potential determined by the unstable
Jacobian, he proved exponential decay of correlations for H\"older
observables under additional hyperbolicity and distortion assumptions
\cite{JakobsonTF17}.  In a related formulation, involving inverse
distortion, bounded initial variation, and a strong-contraction condition,
he obtained exponential mixing for another class of maps with countable
Markov partitions \cite{JakobsonMix17}.  Subsequently,
Jakobson--Simonelli constructed new countable Markov partitions suitable
for thermodynamic formalism and established exponential decay of
correlations and related statistical properties for a related class of
hyperbolic attractors \cite{JS18}.  These works use symbolic thermodynamic
formalism, but their hypotheses and method
are different from the direct operator framework developed here.

A second general route to quantitative statistical properties is provided
by inducing.  Young introduced tower structures that translate hyperbolic
product structure and recurrence estimates into existence of SRB measures
and statistical limit laws \cite{Y98}.  She subsequently related the
tails of return times to rates of mixing \cite{Y99}.  More recently,
Climenhaga--Luzzatto--Pesin constructed SRB measures and Young towers for a
broad class of surface diffeomorphisms \cite{ClLuPe22}.  These results give a
powerful framework for deriving statistical properties from an induced
system.

Our approach here belongs to the functional analytic programme that has become one of the principal methods in modern smooth ergodic theory. The approach goes back, in the one-dimensional setting, to the work of
Lasota--Yorke on piecewise monotone expanding maps
\cite{LY73}. Its decisive extension to
hyperbolic dynamics was initiated by Blank--Keller--Liverani
\cite{BKL02}.  Their work introduced the anisotropic functional-analytic
approach to Anosov maps, in which the Banach space is adapted to the
simultaneous expansion in unstable directions and contraction in stable
directions.  This pioneering construction opened the way for the
development of several anisotropic Banach-space frameworks for smooth
uniformly hyperbolic maps
\cite{Baladi05,GL06,BaTs07,GL08}.  The method was subsequently extended to
piecewise hyperbolic maps, where the dynamics is interrupted by
discontinuity and singularity sets \cite{DeLi08,BaGo09,BG10}, and then to
hyperbolic systems with more general singularity structures, including
the Lorentz gas and broader classes of two-dimensional hyperbolic maps
with singularities \cite{DeZh11,DZ14}.  Surveys and systematic accounts of
this approach can be found in \cite{De18,DKL21}.  The principal strength of
the functional-analytic method is that it translates the hyperbolic
geometry of the system into spectral properties of its transfer
operator: one proves quasi-compactness and the existence of a spectral gap, from which exponential
decay of correlations and a variety of statistical limit theorems
follow. These developments have transformed the study of chaotic dynamical systems into a rich interaction between geometry, analysis, and operator theory.

The purpose of the present paper is to construct anisotropic Banach
spaces adapted to generalized horseshoe maps and to establish a spectral
gap for the associated transfer operator. The elements of these
spaces are distributions whose stable norms are defined by testing their
action against regular functions along stable curves, while a transverse
seminorm compares their action on nearby stable curves.  In this respect,
our norms are closest to the ones in \cite{DeLi08}.  On the resulting spaces, we prove
Lasota--Yorke inequalities and quasi-compactness and 
obtain a spectral gap for the transfer operator. 
Two features prevent a direct application of the existing theory.  First,
the system has countably many branches.  Second, the dynamically defined
stable foliation initially covers only a full-Lebesgue-measure regular
set.  A stable testing family suitable for the transfer operator must
additionally contain stable sides of posts, finite-order singular leaves,
and ordered limiting leaves.  Moreover, uniform regularity is required
both along the leaves and in the transverse direction.

The transverse requirement is substantial and is not an automatic
consequence of uniform hyperbolicity.  For smooth Anosov systems, the
regularity of the invariant splitting is governed by quantitative bunching
conditions, and the corresponding bounds can be sharp
\cite{H94}.  Furthermore, invariant foliations may possess only
low transverse H\"older regularity and need not be transversely Lipschitz
\cite{HW99}.  Thus the strong holonomy estimates obtained below do not follow
from general hyperbolic theory.  They rely on the full-width ordered
geometry of the present model and on uniform distortion estimates for
finite-order stable holonomies.

A classical geometric method for proving the absolute continuity of invariant measures is based on the notion of \emph{transversality}. Roughly speaking, transversality requires that unstable directions associated with different inverse branches remain sufficiently separated, thereby preventing persistent tangencies that obstruct regularity. This approach was developed extensively by Tsujii \cite{T01,T05} and has proved highly successful in the study of expanding and partially hyperbolic systems. More recently, Ban, Oliveira, Liverani and Yang \cite{BrLeObYu25} combined transversality with anisotropic Banach spaces to establish strong spectral properties for transfer operators associated with systems having infinitely many branches.

The notion of \emph{virtual expansion}, introduced by Tsujii \cite{T23}, provides a significant extension of the transversality method. Rather than relying on a global hyperbolic splitting, virtual expansion measures the cumulative expansion of cotangent vectors along all inverse branches in a quantitative and base-point-dependent manner. One of its remarkable features is that the strength of virtual expansion directly controls the Sobolev regularity of the invariant density: stronger virtual expansion implies higher Sobolev regularity and consequently finer statistical properties. This makes virtual expansion particularly well suited for studying systems that are beyond the reach of classical transversality arguments.

In this paper we combine the modern functional analytic approach with Tsujii's theory of virtual expansion to study generalized horseshoe maps. In this 
special case, the ordinary isotropic Sobolev space $H^{\mu}$ ($\mu>0$) already yields quasi-compactness of the transfer operator. 
Under this condition, we  prove that the physical measure is absolutely continuous with respect to Lebesgue measure and that its density belongs to a Sobolev space $H^\mu$ for some $\mu<1/2$. 

\medskip
\hspace{-0.4cm}{\bf Organization of the paper.} The initial Subsection \ref{Int:GHM} is devoted to the definition of generalized horseshoe maps (GHMs) and the conditions imposed on them. The main results of the article are then stated. Section \ref{sec:modula of continuity of stable} is devoted to the study of the geometric and measure-theoretic properties of the stable foliation and to the proof of Theorem \ref{thm:A}. In Section \ref{test:transfer}, we introduce the test spaces and the Banach spaces used throughout the paper. We then study the statistical properties of the transfer operator and, in particular, prove that it has a spectral gap. This yields Theorem \ref{thm:B}. In Section \ref{SecV.E.}, we prove the absolute continuity of the physical measure under the virtually expanding condition, thereby establishing Theorem \ref{thmB}. Finally, Appendix I provides an equivalent characterization of the virtually expanding condition, while Appendix II recalls the necessary preliminaries on Sobolev spaces and establishes a local Lasota--Yorke inequality.
\subsection{Generalized horseshoe maps}\label{Int:GHM}
Let $Q=[0,1]^2$ be the unit square. Fix
$0<\alpha<1$.  The stable and unstable cones are
\[
    K^s_\alpha:=\{(v_1,v_2): |v_1|\le \alpha |v_2|\},
    \qquad
    K^u_\alpha:=\{(v_1,v_2): |v_2|\le \alpha |v_1|\}.
\]
A $C^1$ curve is called a $K^s_\alpha$-curve if it can be written as a
vertical graph
\[
    \{(x(y),y):0\le y\le 1\text{ and }|x'(y)|\le \alpha\}.
\]
Similarly, a $C^1$ curve is called a
$K^u_\alpha$-curve if it can be written as a horizontal graph
\[
    \{(x,y(x)):0\le x\le 1\text{ and } |y'(x)|\le \alpha\}.
\]

We consider a countable collection $\{S_i\}_{i\ge1}$ of closed curvilinear
rectangles contained in $Q$.  Each $S_i$ is a full-height post, that is,
\[
    S_i=\{(x,y)\in Q: x_-^{(i)}(y)\le x\le x_+^{(i)}(y),\ 0\le y\le1\},
\]
where the left and right sides are $K^s_\alpha$-curves. For each $i\ge1$, let
\[
    F_i:\widetilde S_i\longrightarrow F_i(\widetilde S_i)
\]
be a $C^2$ diffeomorphism, where $\widetilde S_i\subset\mathbb R^2$ is an
open domain containing $S_i$.  We require each branch $F_i$ to map $S_i$
diffeomorphically onto a full-width horizontal curvilinear strip, setting
\[
    U_i:=F_i(S_i),
\]
we assume that $U_i$ stretches from the left side of $Q$ to the right side
of $Q$ and has the form
\[
    U_i=\{(x,y)\in Q: y_-^{(i)}(x)\le y\le y_+^{(i)}(x),\ 0\le x\le1\},
\]
where the lower and upper sides are $K^u_\alpha$-curves.  Since each post
is full-height and each image strip is full-width, every transition
$i\to j$ is admissible.  We therefore use the one-sided full countable
shift $\Sigma:=\{\mathbf a=(a_i)_{i\ge1}:a_i\in\N\}.
$
The left shift is $\sigma\mathbf a=(a_2,a_3,\ldots)$. 
For $\mathbf a=(a_1,a_2,\ldots)\in\Sigma$ and $n\ge1$, set
\[
    [\mathbf a]_n:=(a_1,a_2,\ldots,a_n),
\]
The nested vertical subposts are defined inductively by
\[
    S_{[\mathbf a]_1}:=S_{a_1},
    \qquad
    S_{[\mathbf a]_n}
      :=F_{a_1}^{-1}\bigl(S_{[\sigma\mathbf a]_{n-1}}\bigr)
        \cap S_{a_1},
    \qquad n\ge2.
\]
Thus $S_{[\mathbf a]_n}$ consists precisely of the points whose first
n forward branch labels are $a_1,a_2,\ldots,a_n$. These nested
vertical subposts shrink onto the local stable curve associated with
the forward itinerary. The horizontal substrips are defined analogously, but with the finite word read as a backward itinerary. Namely,
\[
    U_{[\mathbf a]_1}:=U_{a_1},
    \qquad
    U_{[\mathbf a]_n}
      :=F_{a_1}\bigl(
          U_{[\sigma\mathbf a]_{n-1}}\cap S_{a_1}
        \bigr),
    \qquad n\ge2.
\]
Thus, in the definition of $U_{[\mathbf a]_n}$, the symbol $a_1$
specifies the present branch, $a_2$ its immediate predecessor, and
so forth. The inner strip, which encodes the earlier part of the
backward itinerary, is constructed first and is then carried forward
by the outermost branch $F_{a_1}$. Consequently, \(U_{[\mathbf a]_{n+1}}\subset U_{[\mathbf a]_n}\subset U_{a_1}\), and these nested full-width
horizontal substrips shrink onto the local unstable curve determined
by the corresponding backward itinerary.

For later use, the branchwise iterate determined by $\mathbf a$ is
\[
    F_{\mathbf a}^n
      :=F_{a_n}\circ\cdots\circ F_{a_2}\circ F_{a_1},
    \qquad n\ge1.
\]
When the symbolic word is fixed, we often write simply $F^n$.

\medskip
\hspace{-.4cm}\textbf{Geometric Conditions.}  We use the Jakobson--Newhouse geometric framework
\cite[104--105]{JN00}.  First,
\[
    \intt(S_i)\cap\intt(S_j)=\varnothing,
    \quad\text{for }i\ne j.
    \tag{G1}
\]
Second,
\[
    \Leb\Big(Q\setminus\bigcup_{i\ge1}S_i\Big)=0.
    \tag{G2}
\]

\medskip
\hspace{-.4cm}\textbf{Hyperbolicity Conditions.}  In the invariant-cone form of
Jakobson--Newhouse, Proposition~2.1
\cite[105--107]{JN00}, there exists $\lambda>1$ such that,
for every $i$,
\[
    DF_i(K^u_\alpha)\subset K^u_\alpha,
    \qquad
    |DF_i v|\ge \lambda|v|
    \quad\text{for every }v\in K^u_\alpha,
    \tag{H1}
\]
and
\[
    DF_i^{-1}(K^s_\alpha)\subset K^s_\alpha,
    \qquad
    |DF_i^{-1}v|\ge \lambda|v|
    \quad\text{for every }v\in K^s_\alpha.
    \tag{H2}
\]
These cone conditions imply that the subposts
$S_{[\mathbf a]_n}$ shrink exponentially in the horizontal direction,
while the substrips $U_{[\mathbf a]_n}$ shrink exponentially in the
vertical direction.

\medskip
\hspace{-.4cm}\textbf{Regularity condition.} Write
$F_i(x,y)=(F_{i,1}(x,y),F_{i,2}(x,y))$, and use $\partial_k F_{i,j}$ and
$\partial_{kl}F_{i,j}$, $j\in\{1,2\}$ and $k,l\in\{x,y\}$, for partial
derivatives.  We define
\[
|D^2F_i(z)|=
\max_{j\in\{1,2\},\, k,l\in\{x,y\}}|\partial_{kl}F_{i,j}(z)|.
\]
We assume the following uniform second-derivative bound:
\[
    B_0:=\sup_{i\ge1}\sup_{z\in S_i}|D^2F_i(z)|<\infty.
\tag{RC1}
\]
This is the uniform \(C^2\) regularity hypothesis used in the piecewise
hyperbolic framework (see \cite[4779--4780]{DeLi08}.

\begin{definition}
A Generalized Horseshoe Map (GHM) on $Q$ is the branchwise piecewise
$C^2$ map generated by the branches $F_i$ above, satisfying $(G1)$--$(G2)$,
$(H1)$--$(H2)$, and $(RC1)$.
It is given on the interior of each post by
\[
    F(z)=F_i(z),
    \qquad z\in\intt(S_i).
\]
On post boundaries, $F$ is interpreted branchwise: a boundary point is
evaluated together with its corresponding branch label.
\label{def:GHM}
\end{definition}

\begin{figure}
	\begin{center}
		\includegraphics[height=4.7cm,width=6cm]{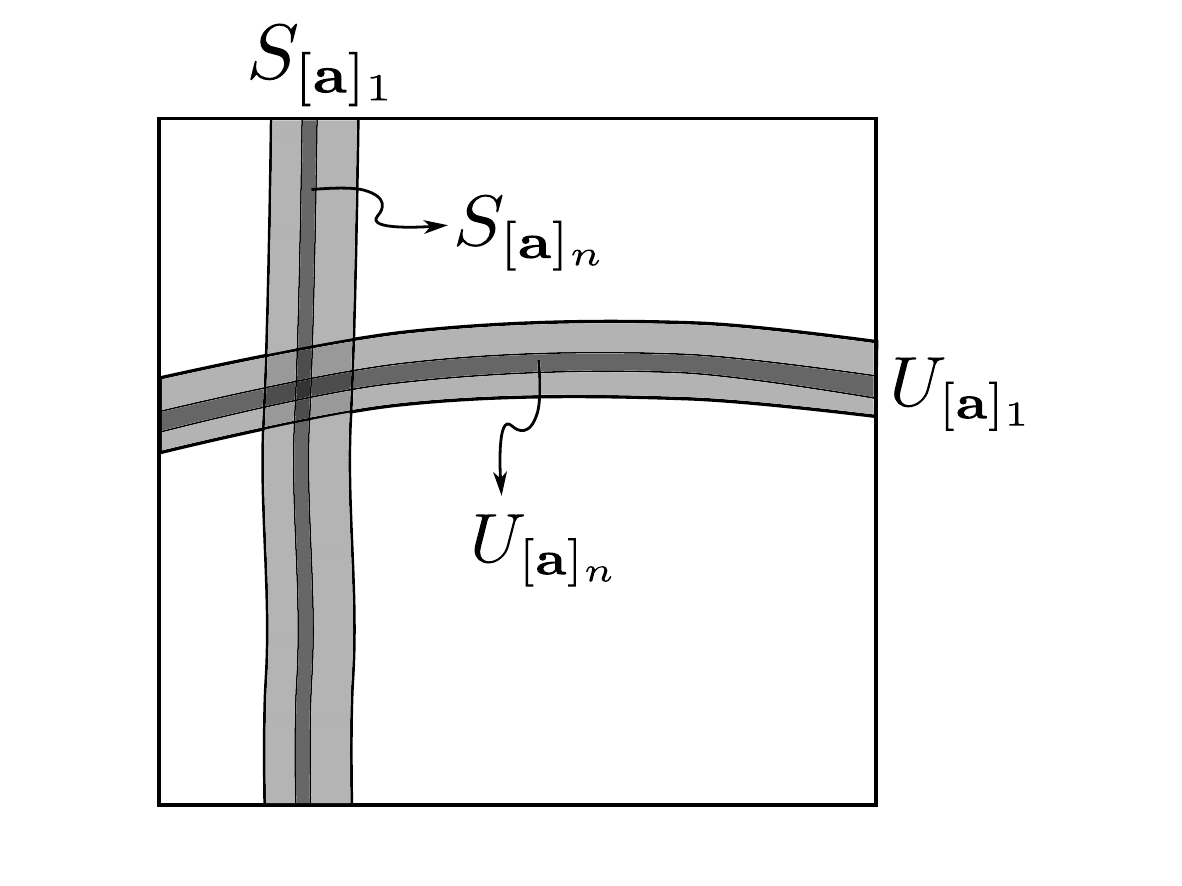}
		\caption{Approximate stable and unstable manifolds\label{StableUnstable} }
	\end{center}
\end{figure}

Under $(H1)$ and $(H2)$, Jakobson--Newhouse,
Proposition~2.2 \cite[pp.~107--108]{JN00}, gives the regular local stable
manifolds as $C^1$ curves tangent to $K^s_\alpha$.
The additional hypothesis (RC1) will be used below to prove the stronger,
uniform $C^{1,1}$ leafwise estimate.  Bounded curvature of admissible unstable
curves is the corresponding key input for the fluctuation estimates of
unstable derivatives and for bounded distortion along unstable curves.
The corresponding bounded-distortion argument follows the composition
estimate of Jakobson--Newhouse, Proposition~8.1
\cite[131--136]{JN00}. For $\mathbf a, \mathbf b\in\Sigma$, define
\[
    W^s_{\loc}(\mathbf a)
    :=\bigcap_{n\ge1}S_{[\mathbf a]_n},
\,\,\,\,\,\,
    W^u_{\loc}(\mathbf b)
    :=\bigcap_{n\ge1}U_{[\mathbf b]_n}.
\]
The stable manifold is a full-height $K^s_\alpha$-curveand and can be written as a graph
$x=x_{\mathbf a}(y)$, whereas, since every branch image is full-width, $W^u_{\loc}(\mathbf b)$ is a
full-width $K^u_\alpha$-curve and can be written as a graph
$y=y_{\mathbf a}(x)$ (see Figure \ref{StableUnstable}). 
For each pair $(\mathbf a,\mathbf b)\in \Sigma\times\Sigma$, the two
curves intersect at a unique point. We denote this point by
\[
    \pi(\mathbf a,\mathbf b)
      :=W^s_{\loc}(\mathbf a)\cap W^u_{\loc}(\mathbf b).
\]
The coded invariant set is then defined by
\[
    \Lambda:=\pi(\Sigma\times\Sigma).
\]
Equivalently, $\Lambda$ is the set of points admitting a full branchwise
symbolic orbit.  By (G1), the posts $S_i$ have pairwise disjoint interiors,
so the forward itinerary of a point is unique away from common boundaries of
posts.  Nothing in $(G1)$--$(G2)$ requires the image strips $U_i$ themselves
to have disjoint interiors: if they overlap, a point may arise as the image
of more than one branch, and so may have several past itineraries.  (Were
the $U_i$ additionally assumed to have pairwise disjoint interiors, the past
itinerary would likewise be unique except on common boundaries of image
strips; we do not impose this here.)  
The topological attractor is defined by
\[
    \Lambda_{\mathrm{top}}
    :=\bigcup_{\mathbf a\in\Sigma}W^u_{\loc}(\mathbf a).
\]
For $z$ on an admissible unstable curve, write
$J^uF^n(z):=|DF^n(z)v|/|v|$ for the expansion factor of $F^n$ at $z$ along a
nonzero tangent vector $v$ to that curve; this is well defined since the
curve is one-dimensional.  Then, for every word $[\mathbf a]_n$ and for
points $z,w$ lying on the same admissible unstable curve inside
$U_{[\mathbf a]_n}$, one obtains
\begin{equation*}\label{BD}
    \frac{J^uF^n(z)}{J^uF^n(w)}\le C,
\end{equation*}
where $C$ is independent of $n$ and of the word.

We recall that an ergodic invariant measure $\mu$ is {\it physical measure} if its basin $B(\mu)$ consists of those points $x\in Q$ for which $1/n\sum_i\delta_{F^i(x)}$ tends to $\mu$ , in the weak topology, has positive Lebesgue measure. In the case that $F$ is non-overlapping, $W^u_{\loc}(\mathbf a)$ defines  a foliation for $\Lambda$. An SRB measure is an invariant probability measure whose conditional measures on local unstable leaves are absolutely continuous with respect to the Riemannian measure induced on those leaves.

\subsection{Statement of the main results}

Our first result concerns the geometry and regularity of the stable
foliation. In addition to the regular stable leaves contained in the
smoothness domains of the map, the natural stable family associated
with a GHM contains finite-order singular leaves and limiting leaves.
Completing the family by including these leaves produces a global
foliation of the square. The following theorem gives the uniform
regularity properties of this completed foliation that will be used
throughout the paper.

\begin{mainthm}
	\label{thm:A}
	The completed stable leaves form a global foliation of the square
	\(Q\). The leaves have uniform \(C^{1,1}\) regularity, the stable
	holonomies between horizontal transversals are uniformly bi-Lipschitz,
	and the stable tangents depend \(1/2\)-Hölder continuously on the
	transverse leaf parameter.
\end{mainthm}

Theorem~\ref{thm:A} provides the geometric structure needed for the
functional-analytic study of the dynamics. In particular, the completed
stable foliation admits a foliated Fubini decomposition with uniformly
controlled conditional densities, while its transverse regularity
allows distributions supported on nearby stable leaves to be compared
uniformly. These properties make it possible to use the completed
stable leaves as a testing family for the construction of weak and
strong anisotropic Banach spaces adapted to the dynamics.

On these spaces, the transfer operator satisfies Lasota--Yorke
inequalities, and the strong space is compactly embedded into the weak
space. This yields quasi-compactness of the transfer operator. To
describe its peripheral spectrum, we consider the quotient obtained by
collapsing each completed stable leaf to a point. The induced
one-dimensional system is a full-branch uniformly expanding map whose
invariant measure is equivalent to Lebesgue measure. Its exactness,
together with the contraction along stable leaves, leads to the
simplicity of the eigenvalue \(1\) and the absence of other peripheral
eigenvalues. The same geometric and measure-theoretic structure also
allows us to identify the physical basin of the invariant measure. We
thus obtain our main statistical result.

\begin{mainthm}
	\label{thm:B}
	Let \(F\) be a {\rm GHM}. There exist anisotropic Banach spaces
	\(
		\mathcal B\hookrightarrow\mathcal B_w
	\)
	such that the transfer operator
	\(
		\mathcal L\colon\mathcal B\longrightarrow\mathcal B
	\)
	is quasi-compact and has a spectral gap. More precisely, \(1\) is a
	simple eigenvalue of \(\mathcal L\) and is its only spectral value on the unit circle. Moreover, \(F\) admits a unique physical invariant probability measure
	\(\mu_0\), whose basin has full planar Lebesgue measure in \(Q\).
	If \(F\) is non-overlapping, then \(\mu_0\) is an SRB measure.
\end{mainthm}

The spectral gap has further statistical consequences. In particular, it implies exponential convergence to equilibrium and exponential decay of correlations for \(C^1\) observables. Moreover, the stability of the strong space under multiplication by \(C^1\)functions allows us to introduce the associated twisted transfer operators. Applying analytic perturbation theory to these operators yields the central limit theorem for \(C^1\)observables. These consequences are stated in Corollary~\ref{thm:statistical_limit_laws}.

To formulate the final result, we assume that $F$ is a finite GHM. We adopt Tsujii’s notion of virtually expanding (see \cite{T23}). 
Let $\mathrm{T}^1S:=S\times \mathbb{S}^1$. For $(x,v)\in \mathrm{T}^1S$ define the projections $\xi_{(x,v)}:\mathbb{R}^2\to \mathbb{R}$ by $\xi_{(x,v)}(w):=\langle v,w\rangle$. Assume that for $x\in S$, $F^{-1}\{x\}:=\{y_1,\ldots,y_n\}$ and put $$\eta(y_i,v):=\xi_{(x,v)}\big(Df(y_i)(B_1(y_i))\big),$$ where $B_1(y_i)$ is the unit ball around $y_i$. Define
$$\beta^\mu(F):=\sup_{(x,v)\in T^1S} b^\mu(x,v),\,\,\text{where}\,\,b^\mu(x,v):=\sum_{i=1}^n \frac{1}{JF(y_i)(\eta(y_i,v))^{\mu}}.$$
Clearly, the sub additivity condition $\beta^\mu(F^{m+n})\leq \beta^\mu(F^n)\,\beta^\mu(F^m)$ holds and one can define
$$B^{\mu}(F):=\lim_{n\to\infty}\big( \beta^{\mu}(F^n)\big)^{1/n}=\inf_{n\in\mathbb{N}}\big( \beta^{\mu}(F^n)\big)^{1/n}.$$

Let $\mathrm{N}(F)$ be the maximum number of preimages of the points in $\Lambda_{\mathrm{top}}$. As before, $\mathrm{N}(F)$ satisfies the sub additivity condition $\mathrm{N}(F^{m+n})\leq \mathrm{N}(F^n)\,\mathrm{N}(F^m)$.  Put 
$$\mathrm{IN}(F):=\lim_{n\to\infty}\big( \mathrm{N}(F)\big)^{1/n}=\inf_{n\in\mathrm{N}}\big(\mathrm{N}(F^n)\big)^{1/n}.$$

We say that an GHM $F:S\to S$ is {\it $\mu$-virtually expanding} if $\sqrt{B^{2\mu}(F)}\,\,\mathrm{IN}(F)<1$. Virtually expanding means $\mu$-virtually expanding for some $\mu>0$ (see Appendix I for the details).
\begin{mainthm} \label{thmB}
Any finite $\mathcal{C}^r$ virtually expanding \textup{GHM}, $r\geq 2$, leaves invariant an absolutely continuous measure whose Radon-Nikodym derivative belongs to a Sobolev space $H^\mu$, for some $\mu<1/2$.
\end{mainthm}
One of the questions that can be raised next is whether a method can be devised in which a higher regularity is achieved? Since a GHM has a Markov structure, meaning that all discontinuity points are mapped to boundary points, we believe that this question is not far-fetched.

In Part II will be devoted to the statistical stability and measures of maximal entropy of GHMs.
\section{Module of Continuity of Stable Foliation}
\label{sec:modula of continuity of stable}

In this section, we study the regularity of the stable foliation, both
along the leaves and in the transverse direction. We first show that
the stable leaves are of class \(C^{1}\), with uniformly Lipschitz derivatives. We
then study the transverse regularity of the foliation, proving
Lipschitz dependence of the leaves and Hölder continuity of their
tangent directions.

\subsection{Full stable foliation}
We first distinguish the dynamically regular stable lamination from
its topological completion. Let
\[
Q_0:=\bigcup_{i\ge1}\operatorname{int}(S_i),
\qquad
Q_n:=Q_0\cap F^{-1}Q_{n-1},
\qquad
\widetilde Q:=\bigcap_{n\ge0}Q_n.
\]
Thus, \(\widetilde Q\) is the full Lebesgue measure set of points whose
entire forward orbits avoid the boundaries of the posts.

We define the regular stable family by
\begin{equation*}
	\mathcal F^s
	:=\left\{
	W^s_{\loc}(\mathbf a):
	\pi(\mathbf a,\mathbf b)\in\widetilde Q
	\text{ for some }\mathbf b\in\Sigma
	\right\}.
	\label{regular_stable_foliation}
\end{equation*}
Since membership in \(\widetilde Q\) is determined by the forward
itinerary, this condition depends only on \(\mathbf a\), while
\(\mathbf b\) merely specifies a point on the corresponding stable
leaf.

The hyperbolicity conditions imply that each member of
\(\mathcal F^s\) is a full-height \(C^1\), \(K_\alpha^s\)-curve and
that distinct leaves are disjoint. Since each leaf is a full-height
graph, the leaves are naturally ordered from left to right. Thus,
\(\mathcal F^s\) forms an ordered family of stable curves. Set
\[
\mathcal W^+
:=\bigcup_{W\in\mathcal F^s}W.
\]
The set \(\mathcal W^+\) is Borel and has full two-dimensional
Lebesgue measure: \(\operatorname{Leb}(Q\setminus\mathcal W^+)=0\). Hence \(\mathcal F^s\) defines a stable lamination of the full-measure
set \(\mathcal W^+\). We next complete this ordered family by filling
the gaps corresponding to irregular itineraries, thereby obtaining a
stable foliation of the whole square \(Q\).

\begin{lemma}
\label{lem:global_straightening}
The ordered family \(\mathcal F^s\) has a unique completion
\[
    \overline{\mathcal F}^s=\{W_t^s:0\le t\le1\}
\]
by full-height \(K_\alpha^s\)-graphs.  There is a continuous function
\(\chi:Q\to[0,1]\) such that
\[
    \chi(t,0)=t,
    \qquad
    |\chi(t,y)-\chi(t,y')|\le\alpha|y-y'|,
\]
and
\[
    W_t^s=\{(\chi(t,y),y):0\le y\le1\}.
\]
The map
\begin{equation}
    \psi:Q\longrightarrow Q,
    \qquad
    \psi(t,y):=(\chi(t,y),y),
\tag{C1}
\label{global_chart}
\end{equation}
is a homeomorphism.
\end{lemma}

\begin{proof}
The proof is simple, but for completeness we include it here. Fix a horizontal transversal \(\Gamma_*=[0,1]\times\{y_*\}\), and let
\(A\subset[0,1]\) be the set of parameters at which some leaf of
\(\mathcal F^s\) meets \(\Gamma_*\).  The set \(A\) is dense. 
For \(t\in A\), write the leaf through \((t,y_*)\) as the graph
\(X_t(y)\), \(0\le y\le1\); then \(X_t(y_*)=t\) and
\(|X_t(y)-X_t(y')|\le\alpha|y-y'|\).  Disjointness orders these graphs:
\(s<t\) implies \(X_s(y)<X_t(y)\) for every \(y\).  Define the missing
graphs by
\begin{equation}
 X(t,y)
 :=\sup_{\substack{s\in A\\s\le t}}X_s(y)
 =\inf_{\substack{s\in A\\s\ge t}}X_s(y).
\tag{3}
\label{completed_stable_graph}
\end{equation}
The \(\sup\) and \(\inf\) agree: a strict
gap at any height would persist on a vertical neighborhood, again leaving
an open rectangle outside \(\mathcal W^+\).  The same argument excludes
jumps of \(t\mapsto X(t,y)\) and gaps at the vertical sides of \(Q\).
Consequently, for every \(y\), \(t\mapsto X(t,y)\) is a continuous,
strictly increasing bijection of \([0,1]\) onto itself.  Moreover,
\[
 |X(t_n,y_n)-X(t,y)|
 \le \alpha|y_n-y|+|X(t_n,y)-X(t,y)|\longrightarrow0,
\]
so \(X\) is jointly continuous.

Consequently \((t,y)\mapsto(X(t,y),y)\) is a continuous bijection of the
square \(Q\) onto itself, hence a homeomorphism, since \(Q\) is compact.
Reparametrizing this map by its value at \(y=0\) -- replacing \(t\) by
\(X(\cdot,0)^{-1}(t)\) -- normalizes it so that \(\chi(t,0)=t\), turning it
into exactly the map \(\psi(t,y)=(\chi(t,y),y)\) of \eqref{global_chart}:
still a homeomorphism of \(Q\), whose vertical fibers
\(W_t^s=\{(\chi(t,y),y):0\le y\le1\}\) are precisely the completed stable
leaves.

If another stable graph were inserted at parameter \(t\), order
preservation would place it between \(X_s\) and \(X_r\) for every
\(s<t<r\), \(s,r\in A\).  Letting \(s\uparrow t\) and \(r\downarrow t\)
and using \eqref{completed_stable_graph} identifies that graph with
\(X(t,\cdot)\).  Thus the completion is unique. 
\end{proof}
\subsection{Regularity of stable leaves}
We write $C^{1,1}([0,1])$ for the space of $C^1$ functions whose first derivative is
Lipschitz.  Thus a $C^{1,1}$ function has a second derivative almost
everywhere. The following lemma shows the uniform leafwise regularity of the complete stable foliation.

\begin{lemma}
\label{lem:completed_leaf_regularity}
Assume $(H1)$, $(H2)$, and $(RC1)$.  There is a constant
$K_s=K_s(\alpha,\lambda,B_0)<\infty$ such that every regular stable leaf,
indexed by \(\mathbf a\in\Sigma\),
\[
    W^s_{\loc}(\mathbf a)
      =\{(x_{\mathbf a}(y),y):0\le y\le1\},
\]
satisfies
\[
    x_{\mathbf a}\in C^{1,1}([0,1]),
    \qquad
    \|x_{\mathbf a}'\|_\infty\le\alpha,
    \qquad
    \Lip(x_{\mathbf a}')\le K_s.
\]
The same conclusions hold for every leaf in the ordered completion
$\overline{\mathcal F}^s$.  Consequently, the straightening function in
Lemma~\ref{lem:global_straightening} satisfies
\begin{equation}
    \chi(t,\cdot)\in C^{1,1}([0,1]),
    \qquad
    |\partial_y\chi(t,y)|\le\alpha,
    \qquad
    \Lip_y(\partial_y\chi(t,\cdot))\le K_s
\tag{C2}
\label{leafwise_C2}
\end{equation}
for every $t\in[0,1]$.
\end{lemma}

\begin{proof}
Fix a regular itinerary
\(\mathbf a=(a_1,a_2,\ldots)\in\Sigma\) and
\(z_0\in W^s_{\loc}(\mathbf a)\).  For each \(n\), intersect the vertical
line through
\[
    z_n=F_{a_n}\circ\cdots\circ F_{a_1}(z_0)
\]
with \(U_{a_n}\), take the component \(\gamma_0\) through \(z_n\), and pull that component
back through \(F_{a_n}^{-1},\ldots,F_{a_1}^{-1}\).  The resulting curve \(\gamma_n\) is a \(C^2\) \(K_\alpha^s\)-graph
\(\{(h_n(y),y):0\le y\le1\}\) contained in
\(S_{[\mathbf a]_n}\).  The horizontal thickness of the \(S_{[\mathbf a]_n}\) is
\(O(\lambda^{-n})\), hence
\[
    \|h_n-x_{\mathbf a}\|_\infty\longrightarrow0.
\label{trial_graph_C0_convergence}
\]
We next obtain a curvature bound independent of \(n\).  For this, we first consider the case $n=1$. That is, we 
consider the graph \(x=g(Y)\)  and its pullback is \(x=h(y)\),
by \(F(x,y)=(F_1(x,y),F_2(x,y))\).  Then
\[
    F_1(h(y),y)=g(Y(y)),
    \qquad Y(y)=F_2(h(y),y).
\]
Put \(p=g'(Y)\), \(v=(h',1)\), and \(\tau=Y'(y)\), with $|\tau|\le\lambda^{-1}$. Then $DF(h',1)=\tau(p,1)$. Put $A:=F_{1x}-pF_{2x}$. The cone conditions give $|A|\ge(1-\alpha^2)\lambda$. Differentiating the graph identity twice gives
\[
    Ah''
      =g''(Y)(Y')^2+pD^2F_2[v,v]-D^2F_1[v,v].
\]
Since \(|h'|,|p|\le\alpha\) and (RC1) gives
\(|D^2F_j[v,v]|\le(1+\alpha)^2B_0\), we obtain
\begin{equation}\label{eq:Hregularity}
    |h''|
       \le a_*+c\,|g''(Y)|,
    \qquad
    a_*:=\frac{(1+\alpha)^3B_0}
                    {(1-\alpha^2)\lambda},
    \qquad
    c:=\frac{\tau^2}{|A|}.
\end{equation}
Now, we apply \eqref{eq:Hregularity} inductively. The essential point is that the transported-curvature coefficients decay
after composition even though a single \(c\) need not be less than one.
Let \(c_j=\tau_j^2/|A_j|\) denotes the coefficient at step \(F_{a_j}\). Then for  \(k\)-step branch composition \(F_{a_1}\circ\ldots\circ F_{a_k}\) its transported-curvature coefficient is therefore
\[
    \frac{\tau_{1,k}^2}{|A_{1,k}|}=\prod_{j=1}^{k}c_j.
\]
Where, \(|\tau_{1,k}|\le\lambda^{-k}\) and $|A_{1,k}|\ge(1-\alpha^2)\lambda^k$. Consequently,
\begin{equation*}
    \prod_{j=0}^{k-1}c_j
       \le\frac{1}{1-\alpha^2}\lambda^{-3k}.
\label{curvature_product_main}
\end{equation*}
Iterating \eqref{eq:Hregularity} yields
\[
    \|h_n''\|_\infty\le K_s,
    \qquad
    K_s:=
    \frac{(1+\alpha)^3B_0}
         {(1-\alpha^2)^2\lambda(1-\lambda^{-3})},
\label{uniform_trial_curvature_main}
\]
uniformly in \(n\) and in the itinerary. Thus \((h_n')\) is uniformly bounded and equi-Lipschitz.  By
Arzel\`a--Ascoli, every subsequence has a uniformly convergent further
subsequence. Hence
\(x_{\mathbf a}\in C^{1,1}\). That is 
\[
\|x_{\mathbf a}'\|_\infty\le\alpha, \quad \textrm{and}\quad
\Lip(x_{\mathbf a}')\le K_s.
\]
%
\end{proof}
\subsection{Lipschitz regularity of stable foliation }
Leafwise regularity alone does not control the transverse parametrization.
Here we prove the uniform transverse bi-Lipschitz regularity of the complete stable foliation. For $r_0,r_1\in[0,1]$, let
\[
    H_{r_0,r_1}:[0,1]\longrightarrow[0,1]
\]
denote the stable holonomy from the horizontal transversal
$\Gamma_{r_0}=[0,1]\times\{r_0\}$ to
$\Gamma_{r_1}=[0,1]\times\{r_1\}$.  Note that with the notation of Lemma~\ref{lem:global_straightening}, 
$H_{r_0,r_1}\bigl(\chi(t,r_0)\bigr)=\chi(t,r_1)$.
\begin{lemma}
\label{lem:transverse_bilip}
Under the assumptions $(H1)$, $(H2)$, and $(RC1)$,
there is a constant $C_H=C_H(\alpha,\lambda,B_0)<\infty$ such that
\begin{equation}
    e^{-C_H|r_1-r_0|}|x-x'|
    \le
    |H_{r_0,r_1}(x)-H_{r_0,r_1}(x')|
    \le
    e^{C_H|r_1-r_0|}|x-x'|
\label{uniform_holonomy_bilip}
\end{equation}
for all $r_0,r_1,x,x'\in[0,1]$.  In particular, putting  $c_\perp=e^{-C_H}$ and $c^\perp=e^{C_H}$, we have
\begin{equation}
    c_\perp|t-s|
    \le|\chi(t,y)-\chi(s,y)|
    \le c^\perp|t-s|,  
\tag{C3}
\label{transverse_C3}
\end{equation}
for all $s,t,y\in[0,1]$.
\end{lemma}

\begin{proof}
Fix a word \(\mathbf a=(a_1,\ldots,a_n)\) and write \(F_{\mathbf a}^n=F_{a_n}\circ\cdots\circ F_{a_1}
\). Let
$I_{\mathbf a}^{(n)}(r)=\{x:(x,r)\in S_{[\mathbf a]_n}\}$ be the horizontal slice of
\(S_{[\mathbf a]_n}\) at height \(r\).  
We first deal with the holonomy map 
\[
    H_{\mathbf a,n}: I_{\mathbf a}^{(n)}(r_0)\to I_{\mathbf a}^{(n)}(r_1).
\]
Write
\[
    DF_i(z)=
    \begin{pmatrix}a(z)&b(z)\\c(z)&d(z)\end{pmatrix}.
\]
The two cone conditions give
\[
    |c(z)|\le\alpha|a(z)|,\qquad
    |b(z)|\le\alpha|a(z)|,\qquad
    |a(z)|\ge\lambda.
\]
For an unstable slope \(u\), define the horizontal expansion and the image slope by
\[
    \Lambda_i(z,u):=a(z)+b(z)u,
    \qquad
    \Psi_i(z,u):=\frac{c(z)+d(z)u}{a(z)+b(z)u}.
\]
Thus \(DF_i(z)(1,u)=\Lambda_i(z,u)(1,\Psi_i(z,u))\), and
\[
    |\Lambda_i(z,u)|\ge(1-\alpha^2)\lambda.
\]
Suppose \(z, z'\) are joined inside \(S_i\) by a \(K^s_{\alpha}\)-curve and \(|u|, |u'|\le\alpha\). The bound (RC1) gives a constant \(C_0=C_0(\alpha,\lambda,B_0)\) such that
\begin{equation}
\begin{aligned}
 |\Psi_i(z,u)-\Psi_i(z',u)|
   &\le C_0d_s(z,z'),\\
 \bigl|\log|\Lambda_i(z,u)|-\log|\Lambda_i(z',u')|\bigr|
   &\le C_0\bigl(d_s(z,z')+|u-u'|\bigr).
\end{aligned}
\label{one_step_holonomy_main}
\end{equation}

We also need contraction in the slope variable over many steps.  Let
\(A=D F^m(z)\) be the derivative of an \(m\)-step branch
composition and write
$A(1,u)=\xi(1,\bar u)
$ with \(|\xi|\ge\lambda^m\).  If
\(\bar w=(p,1)\in K_\alpha^s\) and \(w=A^{-1}\bar w\), then
\(w\in K_\alpha^s\) and \(|w|\ge\lambda^m\). The uniform angle between \(K_\alpha^u\) and \(K_\alpha^s\) implies
\[
|(1,u)\wedge w|\ge (1-\alpha^2)|w|, \quad |(1,\bar u)\wedge\bar w|\le 1+\alpha^2,
\]
and by wedge-product identity
\[
    (\det A)\bigl((1,u)\wedge w\bigr)
      =\xi\bigl((1,\bar u)\wedge\bar w\bigr)
\]
we obtain 
\begin{equation*}\label{eq:det}
    |\det A|
       \le\frac{1+\alpha^2}{1-\alpha^2}\,|\xi|\lambda^{-m}.
\end{equation*}
Since the derivative of the slop map \(\Psi^{(m)}\) associated to \(F^m\) is
\(\det A/\xi^2\), we obtain
\begin{equation}
    |\partial_u\Psi^{(m)}(z,u)|
       \le C_{\rm pr}\lambda^{-2m},
    \qquad
    C_{\rm pr}:=\frac{1+\alpha^2}{1-\alpha^2}.
\label{projective_contraction_main}
\end{equation}
Now let
\[
 z_0=(x,r_0),\qquad
 z_0'=\bigl(H_{\mathbf a,n}(x),r_1\bigr),
 \qquad
 z_k=F_{\mathbf a}^kz_0,\quad z_k'=F_{\mathbf a}^kz_0',
\]
Therefore \[   d_s(z_k,z_k')\le C_s\lambda^{-k}\delta,\] for
$C_s=\sqrt{1+\alpha^2}$ and $\delta:=|r_1-r_0|$. 
Let \(u_k,u_k'\) be the unstable slopes obtained from the two initial
horizontal slices.  Since \(u_0=u_0'=0\), replacing the base points one at
a time, applying \eqref{one_step_holonomy_main}, and propagating each error
with \eqref{projective_contraction_main} gives
\begin{equation*}
\begin{aligned}
 |u_k-u_k'|
 &\le C_0C_{\rm pr}\sum_{j=0}^{k-1}
       \lambda^{-2(k-1-j)}d_s(z_j,z_j')\\
 &\le C_u\lambda^{-k}\delta,
\end{aligned}
\label{paired_slope_main}
\end{equation*}
for \(C_u=C_u(\alpha,\lambda,B_0)\). The horizontal Jacobian of the \(n\)-step  is given by
\begin{equation*}\label{HJacobian}
    J_{\mathbf a,n,r}^x(x)
       =\prod_{k=0}^{n-1}
          |\Lambda_{a_{k+1}}(z_k,u_k)|.
\end{equation*}
Using the second estimate in \eqref{one_step_holonomy_main}, followed by
the position and slope estimates above, yields
\begin{equation*}
 \Big|
   \log\frac{J_{\mathbf a,n,r_0}^x(x)}
   {J_{\mathbf a,n,r_1}^x
      \big( H_{\mathbf a,n}(x)\big)}
 \Big|
 \le C_H\delta,
\end{equation*}
where \(C_H=C_H(\alpha,\lambda,B_0)\) is is uniform in \(\mathbf a, n, r_0\) and \(r_1\). Consequently,
\begin{equation}\label{Jc:hol}
    1\simeq \frac{\big|F_{\mathbf a}^n(I_{\mathbf a}^{(n)}(r_1))\big|}
    {\big|F_{\mathbf a}^n(I_{\mathbf a}^{(n)}(r_0))\big|}
    \leq C_H\delta\,\,\frac{\big|I_{\mathbf a}^{(n)}(r_1)\big|}{\big|I_{\mathbf a}^{(n)}(r_0)\big|}
    = C_H\delta\,\,\frac{\big|H_{\mathbf a,n}(I_{\mathbf a}^{(n)}(r_1))\big|}{\big|I_{\mathbf a}^{(n)}(r_0)\big|}
\end{equation}
Now, \eqref{Jc:hol} implies that 
\begin{equation}
    e^{-C_H\delta}\big |I_{\mathbf a}^{(n)}(r_0)\big|
       \le
          \big |H_{\mathbf a,n}(I_{\mathbf a}^{(n)}(r_0))\big |      \le e^{C_H\delta}\big|I_{\mathbf a}^{(n)}(r_0)\big| .
\label{finite_holonomy_main}
\end{equation}
It remains to pass to the completed holonomy. The unstable expansion in (H1) gives \(|I_{\mathbf a}^{(n)} (r)| \le C_\alpha \lambda^{-n}\) uniformly in \(\mathbf a\) and \(r\), and hence the cylinder intervals \(I_{\mathbf a}^{(n)} (r)\) generate the Borel
\(\sigma\)-algebra.  Both the completed holonomy and the finite-cylinder
holonomy send \(I_{\mathbf a}^{(n)} (r_0)\) to \(I_{\mathbf a}^{(n)} (r_1)\), up to endpoints.  Thus
\eqref{finite_holonomy_main} gives
\[
 e^{-C_H\delta}|E|
 \le|H_{r_0,r_1}(E)|
 \le e^{C_H\delta}|E|
\]
for every
Borel set \(E\). In particular, taking \(E=[x,x']\) gives
\[
 e^{-C_H|r_1-r_0|}|x-x'|
 \le|H_{r_0,r_1}(x)-H_{r_0,r_1}(x')|
 \le e^{C_H|r_1-r_0|}|x-x'|.
\]
This is \eqref{uniform_holonomy_bilip}; the height-independent constants
\(e^{\pm C_H}\) follow because \(|r_1-r_0|\le1\).  Finally,
\(H_{0,y}(t)=\chi(t,y)\), so \eqref{transverse_C3} follows.
\end{proof}
\begin{remark}
Jakobson--Newhouse in \cite{JN00} proved the $C^1$ regularity of stable leaves and absolute continuity of the holonomy. In our results, the role of condition $(RC1)$ is crucial to get more regularity in Lemmas \ref{lem:completed_leaf_regularity} and \ref{lem:transverse_bilip}.
\end{remark}
The final regularity statement follows from (C2) and (C3) by interpolation. We use the norm
\[
    \|f\|_{C^1([0,1])}:=\|f\|_\infty+\|f'\|_\infty.
\]

\begin{lemma}
\label{lem:transverse_tangent_holder}
For all $s,t\in[0,1]$,
\begin{equation}
    \|\chi(t,\cdot)-\chi(s,\cdot)\|_{C^1([0,1])}
       \le C_{1/2}|t-s|^{1/2},
\tag{C4}
\label{transverse_C4}
\end{equation}
where $ C_{1/2}
      :=9c_\perp+2\sqrt{2K_sc_\perp}<\infty
$
\end{lemma}

\begin{proof}
Put
\[
 h(y):=\chi(t,y)-\chi(s,y),
 \qquad \delta:=|t-s|.
\]
By (C3), \(\|h\|_\infty\le c_\perp\delta\); and since \(\partial_y\chi(t,
\cdot)\) and \(\partial_y\chi(s,\cdot)\) each have Lipschitz constant at
most \(K_s\) by (C2), so does their difference \(h'\), which gives
\(\Lip(h')\le2K_s\).

We now bound \(\|h'\|_\infty\). Fix \(y\in[0,1]\) and
\(0<\ell\le\frac12\); since \(y\) is within \(\frac12\) of one of the two
endpoints of \([0,1]\), at least one of the two intervals
\([y,y+\ell]\) or \([y-\ell,y]\) is contained in \([0,1]\) -- say it is
\([y,y+\ell]\), the other case being identical with a sign change.
Writing
\[
    h(y+\ell)-h(y)=\int_0^\ell h'(y+r)\,dr
    =\ell\,h'(y)+\int_0^\ell\bigl(h'(y+r)-h'(y)\bigr)\,dr,
\]
and using \(|h(y+\ell)-h(y)|\le2\|h\|_\infty\) together with
\(|h'(y+r)-h'(y)|\le\Lip(h')\,r\) for the integrand, gives
\[
    |h'(y)|
    \le\frac{2\|h\|_\infty}{\ell}
       +\frac{\Lip(h')}{2}\ell,
    \qquad\text{for every } 0<\ell\le\tfrac12.
\]
This holds for every admissible \(\ell\), so we may choose \(\ell\) to
minimize the right side.  If \(\Lip(h')=0\), then \(h'\) is constant and
\(|h'|=|h(1)-h(0)|\le2\|h\|_\infty\).  Otherwise, its
minimizer is \(\ell^*=2\sqrt{\|h\|_\infty/\Lip(h')}\). If
\(\ell^*\le\frac12\), take \(\ell=\ell^*\) and if instead \(\ell^*>\frac12\), i.e.\
\(\Lip(h')<16\|h\|_\infty\), take \(\ell=\frac12\).
In either case,
\[
    \|h'\|_\infty
    \le2\sqrt{\|h\|_\infty\Lip(h')}+8\|h\|_\infty.
\]
Substituting the bounds \(\|h\|_\infty\le c_\perp\delta\) and
\(\Lip(h')\le2K_s\) obtained above gives
\[
    \|h'\|_\infty
    \le2\sqrt{2K_sc_\perp\delta}+8c_\perp\delta
    \le (2\sqrt{2K_sc_\perp}+8c_\perp)\,\delta^{1/2}.
\]
Adding the two displayed bounds,
\[
    \|h\|_{C^1([0,1])}
    =\|h\|_\infty+\|h'\|_\infty
    \le\Bigl(9c_\perp+2\sqrt{2K_sc_\perp}\Bigr)\delta^{1/2}
    =C_{1/2}\delta^{1/2}.
\]
This is exactly \eqref{transverse_C4}.
\end{proof}
\begin{remark}
Properties (C1)--(C4) use the uniform $C^2$ bound (RC1).  The
exponent $1/2$ in (C4) is the interpolation exponent
provided by a uniform leafwise $C^{1,1}$ bound and transverse $C^0$
Lipschitz regularity.  Stronger transverse conclusions may require an
additional bunching hypothesis.
\end{remark}
\subsection{Disintegration along stable leaves}
The analogue of ordinary Fubini in the straight vertical foliation requires
the conditional densities of planar Lebesgue measure along the curved stable
leaves.  Jakobson--Newhouse prove absolute continuity of this stable
foliation in Proposition~10.1 \cite[pp.~140--145]{JN00}; the following lemma
records the stronger quantitative form needed below.

\begin{lemma}
\label{lem:foliated_Fubini}
There are constants \(0<c_{\rm dis}\le C_{\rm dis}\) and a bounded
Borel function $\rho:Q\to(0,\infty)$ such that, upon writing
$\rho_t:=\rho|_{W_t^s}$,
\begin{equation}
	\sup_{t\in[0,1]}|\rho_t|_{C^1(W_t^s)}
	\le C_{\rm dis},
	\qquad
	\inf_{t\in[0,1]}\inf_{z\in W_t^s}\rho_t(z)
	\ge c_{\rm dis}
	\label{conditional_density_bounds}
\end{equation}
such that for every $u\in L^1(Q,m)$,
\begin{equation*}
	\int_Q u\,dm
	=\int_0^1\int_{W_t^s}
	u\,\rho_t\,dm_t\,dt.
	\label{foliated_Fubini}
\end{equation*}
where $m_t$ is the Lebesgue measure on $W_t^s$. 
\end{lemma}

\begin{proof}
It is enough to note that
\begin{equation*}
	\rho_t(\psi(t,y))
	=\frac{\partial_t\chi(t,y)}
	{\sqrt{1+|\partial_y\chi(t,y)|^2}},
\label{conditional_density_formula}
\end{equation*}
Now, in view of Lemma \ref{lem:completed_leaf_regularity} and \ref{lem:transverse_bilip}, a direct calculation give the desired constants $c_{\textrm{dis}}$ and $C_{\textrm{dis}}$ in \eqref{conditional_density_bounds}.
\end{proof}
\section{Test Space and Transfer Operator}\label{test:transfer}
For \(n\ge1\), let \(\mathcal S_n\) be the family of left and right \(K_\alpha^s\)-curves bounding each \(S_{[\mathbf a]_n}\),
\(\mathbf a\in\Sigma\), and put
\[
    \mathcal S:=\bigcup_{n\ge1}\mathcal S_n.
\]
These finite-order singular curves and their ordered limiting curves occur in
the completion \(\overline{\mathcal F}^s\). We henceforth take
$
    \mathscr W^s:=\overline{\mathcal F}^s
$
as the stable testing family.  Lemma~\ref{lem:completed_leaf_regularity}
shows that all regular, finite-order singular, and limiting leaves in this
family are uniformly $C^{1,1}$.  Lemma~\ref{lem:transverse_bilip} moreover
gives uniform bi-Lipschitz stable holonomy between horizontal transversals.
Finally, Lemma~\ref{lem:transverse_tangent_holder} gives uniform
$1/2$-H\"older regularity of stable tangents with respect to the transverse
leaf parameter. The family \(\mathscr W^s\) is interpreted branchwise on singular curves.  If
\(W\in\mathscr W^s\) carries the branch label \(i\), so that \(W\subset S_i\),
then there is \(W^+\in\mathscr W^s\) such that
$F_i(W)\subset W^+
$. 
For a regular stable leaf this is the usual invariance
$
    F_{a_1}\bigl(W^s_{\loc}(\mathbf a)\bigr)
    \subset W^s_{\loc}(\sigma\mathbf a).
$ 
Similarly, branchwise preimages decompose into stable testing curves carrying
their natural labels.

\subsection{Leafwise test spaces}


For $W\in\mathscr W^s$, let $\Lip(W)$ denote the space of Lipschitz functions
on $W$, with norm
\[
    |\varphi|_{\Lip(W)}:=|\varphi|_{C^0(W)}+\Lip_W(\varphi),\qquad |\varphi|_{C^1(\mathscr W^s)}:=\sup_{W\in\mathscr W^s}|\varphi|_{\Lip(W)}
\]
For $0<\gamma<1$, define
\[
    |\varphi|_{C^\gamma(W)}
    :=|\varphi|_{C^0(W)}+H^\gamma_W(\varphi), \qquad |\varphi|_{C^\gamma(\mathscr W^s)}:=\sup_{W\in\mathscr W^s}|\varphi|_{C^\gamma(W)}
\]
where
\[
    H^\gamma_W(\varphi)
    :=\sup_{z\ne z'\in W}
    \frac{|\varphi(z)-\varphi(z')|}{d_W(z,z')^\gamma}.
\]
Define $C^1(\mathscr W^s)$ to be the space of \emph{leafwise Lipschitz} functions.  More precisely,
\[
 C^1(\mathscr W^s):=
 \left\{\varphi:Q\to\mathbb C:
 \varphi\text{ is bounded and Borel measurable, and }
 |\varphi|_{C^1(\mathscr W^s)}<\infty\right\},
\]
Note that this set may be smaller than the class of all
bounded Borel functions having finite \(|\cdot|_{C^\gamma(\mathscr W^s)}\)-norm. Put 
\[    
    C^\gamma(\mathscr W^s)
    :=\overline{C^1(\mathscr W^s)}^{\,|\cdot|_{C^\gamma(\mathscr W^s)}}
\]

The forward stability and contraction of
\(\mathscr W^s\) imply that the branchwise Koopman action preserves these
completed test spaces. To verify this 
begin with \(\varphi\in C^\gamma(\mathscr W^s)\).  
If \(W\in\mathscr W^s\) is labelled by
branch \(i\) and \(F_i(W)\subset W^+\), then for \(z,z'\in W\),
\[
    |\varphi(F_i z)-\varphi(F_i z')|
    \le H^\gamma_{W^+}(\varphi)d_{W^+}(F_i z,F_i z')^\gamma\\
    \le \lambda^{-\gamma}H^\gamma_{W^+}(\varphi)d_W(z,z')^\gamma,
\]
and, so
\begin{equation}\label{koopman}
    |\varphi\circ F|_{C^\gamma(\mathscr W^s)}
    \le|\varphi|_{C^\gamma(\mathscr W^s)}.
\end{equation}
Moreover, stable contraction shows directly that
\(\varphi\circ F\in C^1(\mathscr W^s)\), with branch labels retained on
singular curves.  
Now let
\(\varphi_n\in C^1(\mathscr W^s)\) converge to
\(\varphi\in C^\gamma(\mathscr W^s)\) in the \(C^\gamma\)-norm.  Applying
\eqref{koopman} to \(\varphi_n-\varphi_m\) shows that
\((\varphi_n\circ F)\) is Cauchy in \(C^\gamma(\mathscr W^s)\).  Its limit is
independent of the approximating sequence and is, by uniform convergence,
the branchwise composition \(\varphi\circ F\). Thus the Koopman action
extends continuously to \(C^\gamma(\mathscr W^s)\). 
The backward stability of \(\mathscr W^s\) will be used when testing iterates
of the transfer operator: preimages of a stable testing curve decompose into
stable testing curves.
\subsection{The transfer operator}


For an integrable density $h$, the transfer operator associated with $F$ is
given almost everywhere by
\[
    (\mathcal L h)(z)
    =\sum_{i:\, z\in U_i}
    \frac{h(F_i^{-1}z)}{|\det DF_i(F_i^{-1}z)|}.
\]
The sum is taken over all full-width branches whose images contain $z$. 
Indeed, Tonelli's theorem and branchwise change of variables give
\[
 \int_Q\sum_{i:\,z\in U_i}
 \frac{|h(F_i^{-1}z)|}{|\det DF_i(F_i^{-1}z)|}\,dm(z)
 =\sum_{i\ge1}\int_{S_i}|h(w)|\,dm(w)\\
 =\int_Q|h|\,dm.
\]
Thus the series of absolute values is finite for almost every $z$, and
$\|\mathcal Lh\|_{L^1}\le\|h\|_{L^1}$. Equivalently, and more robustly for distributions, we define
\[
    \mathcal L h(\varphi)=h(\varphi\circ F),
    \qquad
    \varphi\in C^\gamma(\mathscr W^s).
\]
Because of \eqref{koopman}, the right-hand side is well defined whenever $h$ acts continuously on
$C^\gamma(\mathscr W^s)$.  The same formula is read with
$\varphi\in C^1(\mathscr W^s)$. 
 For integrable densities $h$, the two definitions agree.  Moreover, the mass
functional is preserved, that is
\[
    (\mathcal L h)(1)=h(1).
\]
When $h$ is a density with respect to Lebesgue measure, this reads
\[
    \int_Q \mathcal L h\,dm=\int_Q h\,dm.
\]

\subsection{The Banach spaces}
We take the seed space to be
\begin{equation}\label{Def:D}
    \mathcal D:=C^1(Q),
\end{equation}
where $C^1(Q)$ means the space of functions on $Q$ which extend to $C^1$
functions on an open neighborhood of $Q$.  The anisotropic spaces below are
obtained by completing this smooth seed space in the corresponding norms.
Fix, once and for all,
\[
    0<\gamma<1,\qquad
    0<\beta<\rho:=\frac{1-\gamma}{2},\qquad
\]
The restriction $\beta<\rho$ is the transverse bunching required by the
strong unstable estimate.  The exponent $1/2$ comes from (C4), while the
factor $1-\gamma$ is the loss in passing from a transverse $C^0$ estimate to
a $C^\gamma$ estimate by interpolation.

\medskip
\hspace{-0.4cm}\textbf{The weak norm.}  For $h\in\mathcal D$, define
\[
    |h|_w
    =\sup_{W\in\mathscr W^s}
    \sup_{\substack{\varphi\in\Lip(W)\\
        |\varphi|_{\Lip(W)}\le1}}
    \left|
    \int_W h\,\varphi\,dm
    \right|.
\]
This norm tests $h$ leafwise against Lipschitz functions on stable curves.  The
weak space is then defined by completion:
\[
    \mathcal B_w:=\overline{\mathcal D}^{\,|\cdot|_w}.
\]

\medskip
\hspace{-0.6cm}
\textbf{The strong stable norm.}  Fix $0<\gamma<1$.  The strong stable norm is
\[
    \|h\|_s
    =\sup_{W\in\mathscr W^s}
    \sup_{\substack{\varphi\in C^\gamma(W)\\
        |\varphi|_{C^\gamma(W)}\le1}}
    \left|
    \int_W h\,\varphi\,dm_W
    \right|.
\]
Since the test functions are less regular than Lipschitz functions, this norm
is stronger than the weak norm.

\medskip
\hspace{-0.6cm}
\textbf{The strong unstable norm.}  This norm detects regularity in the
unstable direction by comparing integrals over nearby stable curves. Let \(\psi\) be the homeomorphism constructed in Lemma
\ref{lem:global_straightening}.  For
\[
W_i:=W_{t_i}^s=\psi_{t_i}([0,1]),
\qquad i=1,2,
\]
define
\[
d_{\mathscr W^s}(W_1,W_2):=|t_1-t_2|.
\]
Because \(\chi(t,0)=t\), this is a genuine metric on the completed leaves.
By \eqref{transverse_C3},
\[
|t_1-t_2|
\le \sup_{y\in[0,1]}|\chi(t_1,y)-\chi(t_2,y)|
\le c^\perp|t_1-t_2|.
\]
Thus the chart metric is uniformly equivalent to the $C^0$ graph distance. For \(W_t^s=\psi_t([0,1])\) and \(\varphi\in\Lip(W_t^s)\), put
$
\widetilde\varphi:=\varphi\circ\psi_t
$
and define
\[
d_{C^\gamma}(\varphi_1,\varphi_2)
:=|\widetilde\varphi_1-\widetilde\varphi_2|_{C^\gamma([0,1])}.
\]
Moreover, if
\[
dm_{W_t^s}=g(t,y)\,dy,
\qquad
g(t,y)=\sqrt{1+|\partial_y\chi(t,y)|^2},
\]
then by (C4), 
\begin{equation}\label{H1/2}
\|g(t_1,\cdot)-g(t_2,\cdot)\|_\infty
\le C_{1/2}|t_1-t_2|^{1/2}.
\end{equation}
Thus both the graph parametrizations and their arclength densities have the
transverse control needed when integrals over nearby stable leaves are
compared. Define
\[
	\|h\|_u
	:=\sup_{\substack{W_i\in\mathscr W^s,\ i=1,2\\
			0<d_{\mathscr W^s}(W_1,W_2)\le1}}
	\sup_{\substack{\varphi_i\in\Lip(W_i)\\
			|\varphi_i|_{\Lip(W_i)}\le1\\
			d_{C^\gamma}(\varphi_1,\varphi_2)
			\le d_{\mathscr W^s}(W_1,W_2)}}\hspace{-0.4cm}
	\frac{1}{d_{\mathscr W^s}(W_1,W_2)^\beta}
	\left|
	\int_{W_1}h\varphi_1\,dm_{W_1}
	-\int_{W_2}h\varphi_2\,dm_{W_2}
	\right|.
\]
This is finite on the seed space.  Indeed, after pulling both leaf integrals
back to $[0,1]$, the $C^0$ closeness of the paired test functions and (C3)
give errors of order $d_{\mathscr W^s}(W_1,W_2)$, while
\eqref{H1/2} gives an error of order
$d_{\mathscr W^s}(W_1,W_2)^{1/2}$.  Hence, for $h\in\mathcal D$,
\[
\left|
\int_{W_1}h\varphi_1\,dm_{W_1}
-\int_{W_2}h\varphi_2\,dm_{W_2}
\right|
\le C\|h\|_{C^1(Q)}
d_{\mathscr W^s}(W_1,W_2)^{1/2}.
\]
The choice $\beta<1/2$ therefore ensures $\|h\|_u<\infty$. The strong anisotropic norm is
\[
\|h\|_{\mathcal B}:=\|h\|_s+\|h\|_u.
\]
The strong space is
\[
\mathcal B:=\overline{\mathcal D}^{\,\|\cdot\|_{\mathcal B}}.
\]

\subsection{Summability bound.} The countable-branch estimates below require a uniform summability bound for
the geometric weights produced by branchwise change of variables.  We record
that bound here before applying it to the transfer operator. 

Fix \(n\ge1\). For each word
\(\mathbf a=(a_1,\ldots,a_n)\in\mathbb N^n\), let
\(S_{[\mathbf a]_n}\) denote the corresponding vertical subpost and recall that
\(
    I_{\mathbf a}^{(n)}(0)
      :=\{x\in[0,1]:(x,0)\in S_{[\mathbf a]_n}\}
\).
For \(W\in\mathscr W^s\), let \(\mathcal G_n(W)\) denote the collection
of its branchwise \(n\)-step preimages:
\[
    \mathcal G_n(W)
      :=\left\{
          W_{\mathbf a}
          :=
          (F_{\mathbf a}^n)^{-1}
          \bigl(
             W\cap F_{\mathbf a}^n(S_{[\mathbf a]_n})
          \bigr)
          :
          \mathbf a\in\mathbb N^n
        \right\}.
\]
The full-width crossing property implies that each
\(W_{\mathbf a}\) belongs to \(\mathscr W^s\). If
\(J_{W_{\mathbf a}}F_{\mathbf a}^n\) denotes the arclength Jacobian of
\(F_{\mathbf a}^n\) along \(W_{\mathbf a}\), we define the corresponding
inverse transverse Jacobian  by
\begin{equation}
    \mathfrak j_{\mathbf a,W}^{(n)}(z)
    :=
    \frac{J_{W_{\mathbf a}}F_{\mathbf a}^n(z)}
         {|\det DF_{\mathbf a}^n(z)|},
    \qquad z\in W_{\mathbf a}.
    \label{inverse_branch_weight}
\end{equation}

\begin{lemma}
	\label{lem:inverse_branch_weight_summability}
	There is a constant $C_J=C_J(\alpha,\lambda,B_0)\ge1$ such that
	\begin{equation}
		|\mathfrak j_{\mathbf a,W}^{(n)}|_{\Lip(W_{\mathbf a})}
		\le C_J|I_{\mathbf a}^{(n)}(0)|
		\label{individual_inverse_weight_bound}
	\end{equation}
	for every $n$, $W$, and $W_{\mathbf a}\in\mathcal G_n(W)$.  Consequently,
	\begin{equation}
		\sup_{n\ge1}\sup_{W\in\mathscr W^s}
		\sum_{\mathbf a\in\mathbb{N}^n}
		|\mathfrak j_{\mathbf a,W}^{(n)}|_{\Lip(W_{\mathbf a})}
		\le C_J.
		\label{inverse_weight_summability}
	\end{equation}
\end{lemma}

\begin{proof}
	To begin the proof, first note that
    \begin{equation}
		\left|
		\log\mathfrak j_{\mathbf a,W}^{(n)}(z)
		-\log\mathfrak j_{\mathbf a,W}^{(n)}(z')
		\right|
		\le C_s^*d_{W_{\mathbf a}}(z,z').
		\label{stable_weight_distortion}
	\end{equation}
    For any $z,z'$ belong to the same stable preimage $W_{\mathbf a}$. This follows by the same calculation used in Lemma~\ref{lem:transverse_bilip} to \(J_{\mathbf a,n}^x\). 
    
    The strategy is to first compare inverse transverse Jacobian with the inverse horizontal expansion and then explore the connection between the inverse horizontal expansion and the lengths of the intervals. Fix
	$z\in W_{\mathbf a}$, set $A_z=DF_{\mathbf a}^n(z)$, and let $\tau(z)$ be
	the unit tangent to $W_{\mathbf a}$ at $z$.  Cone invariance gives \(Ae_x\in K_\alpha^u\) and \(A\tau^s(z)\in K_\alpha^s\). For the horizontal Jacobian $J_{\mathbf a,n}^x(z):=J_{\mathbf a,n,0}^x(z)$ given by \eqref{HJacobian}, applying \eqref{eq:det}, we get
		\begin{equation}
		\frac{C_{\rm pr}^{-1}}{J_{\mathbf a,n}^x(z)}\le \mathfrak j_{\mathbf a,W}^{(n)}(z)=\frac{|A_z\tau(z)||e_x\wedge\tau(z)|}{|A_ze_x\wedge A_z\tau(z)|}=\frac{|e_x\wedge\tau(z)||\pi_x u_n(z)|}{J_{\mathbf a,n}^x(z)|u_n(z)\wedge v_n(z)|}
		\le \frac{C_{\rm pr}}{J_{\mathbf a,n}^x(z)},
		\label{weight_unstable_comparison}
	\end{equation}
     where \( u_n(z)=\frac{A_ze_x}{|A_ze_x|}, v_n(z)=\frac{A_z\tau(z)}{|A_z\tau(z)|}\) and $\pi_x$ is the projection along the vertical lines. The last inequality holds since
     \[J_{\mathbf a,n}^x(z)=\left|\partial_x(\pi_xF_{\mathbf a}^n)(z)\right|=|\pi_xA_ze_x|=|A_ze_x||\pi_xu(z)|.\]
	
	.
	Now, we compare horizontal Jacobian by the lengths. By \eqref{BD}, the horizontal Jacobian has uniform distortion on a horizontal slice, that is
	\begin{equation}
		\sup_{x,x'\in I_{\mathbf a}^{(n)}(0)}
		\left|
		\frac{J_{\mathbf a,n}^x(x,0)}
		{J_{\mathbf a,n}^x(x',0)}
		\right|
		\le C_u^*.
		\label{horizontal_cylinder_distortion}
	\end{equation}
Note that the map \(\pi_xF_{\mathbf a}^n(x,0)\) is a diffeomorphism from $I_{\mathbf a}^{(n)}(0)$ onto $[0,1]$.  Hence
	\[
	1=\int_{I_{\mathbf a}^{(n)}(0)}
	J_{\mathbf a,n}^x(x,0)\,dx.
	\]
	Together with \eqref{horizontal_cylinder_distortion}, this gives
	\begin{equation}
		\sup_{x\in I_{\mathbf a}^{(n)}(0)}
		\frac1{J_{\mathbf a,n}^x(x,0)}
		\le {C_u^*}|I_{\mathbf a}^{(n)}(0)|.
		\label{inverse_horizontal_width_bound}
	\end{equation}
    Equations
	\eqref{stable_weight_distortion}, \eqref{weight_unstable_comparison} and
	\eqref{inverse_horizontal_width_bound} show that
	\[
	\|\mathfrak j_{\mathbf a,W}^{(n)}\|_\infty
	\le C|I_{\mathbf a}^{(n)}(0)|.
	\]
	Moreover, \eqref{stable_weight_distortion} and the elementary inequality
	$|e^r-1|\le e^{|r|}|r|$ imply
	\[
	\Lip_{W_{\mathbf a}}(\mathfrak j_{\mathbf a,W}^{(n)})
	\le C\|\mathfrak j_{\mathbf a,W}^{(n)}\|_\infty.
	\]
	This proves \eqref{individual_inverse_weight_bound}.  Finally, the intervals
	$I_{\mathbf a}^{(n)}(0)$ associated with distinct words have disjoint
	interiors.  Therefore
	\[
	\sum_{\mathbf a\in\mathbb{N}^n}|I_{\mathbf a}^{(n)}(0)|\le1,
	\]
	and \eqref{inverse_weight_summability} follows.
\end{proof}

\subsection{Lasota–Yorke estimates. }

We now organize the transfer operator branch by branch. For
\(h\in\mathcal D\), \(n\ge1\), and
\(\mathbf a\in\mathbb N^n\), define the \(n\)-step branch operator by
\[
\mathcal L_{\mathbf a}^n h
:=
\mathbf1_{F_{\mathbf a}^n(S_{[\mathbf a]_n})}
\frac{h\circ(F_{\mathbf a}^n)^{-1}}
{|\det DF_{\mathbf a}^n|\circ(F_{\mathbf a}^n)^{-1}}.
\]
Accordingly, the \(n\)-step transfer operator admits the branchwise
decomposition
\[
\mathcal L^n h
=
\sum_{\mathbf a\in\mathbb N^n}\mathcal L_{\mathbf a}^n h.
\]
We first establish the weak and strong norm estimates for
\(h\in\mathcal D\). These estimates will subsequently allow us to
extend the transfer operator to the completed spaces
\(\mathcal B_w\) and \(\mathcal B\).

\begin{lemma}
	\label{lem:weak_norm_estimate}
	There is a constant $C_w$, independent of $n$, such that the transfer
	operator satisfies
	\begin{equation}
		|\mathcal L^n h|_w\le C_w|h|_w,
		\qquad h\in\mathcal D,\quad n\ge0.
		\label{weak_norm_estimate}
	\end{equation}
\end{lemma}

\begin{proof}
	Fix $h\in\mathcal D$, $W\in\mathscr W^s$ and
	$|\varphi|_{\Lip(W)}\le1$.  Changing variables along each branchwise inverse image
	$W_{\mathbf a}\in\mathcal G_n(W)$ gives
	\begin{equation}
		\int_W\mathcal L^nh\,\varphi\,dm
		=\sum_{\mathbf a\in\mathbb{N}^n}
		\int_{W_{\mathbf a}}h\,
		(\varphi\circ F_{\mathbf a}^n)
		\mathfrak j_{\mathbf a,W}^{(n)}\,dm,
		\label{weak_change_variables}
	\end{equation}
	Stable contraction implies
	$|\varphi\circ F_{\mathbf a}^n|_{\Lip(W_{\mathbf a})}\le1$.
	Therefore by \eqref{inverse_weight_summability},
	\[
	\left|\int_W\mathcal L^nh\,\varphi\,dm\right|
	\le |h|_w
	\sum_{\mathbf a\in\mathbb{N}^n}
	|\mathfrak j_{\mathbf a,W}^{(n)}|_{\Lip(W_{\mathbf a})}\le C_J|h|_w
	.
	\]
	Taking the suprema over $W$ and $\varphi$ proves
	\eqref{weak_norm_estimate}, with $C_w=C_J$.
\end{proof}

\begin{lemma}
	\label{lem:strong_stable_estimate}
	There is a constant $C_s$, independent of $n$, such that
	\begin{equation}
		\|\mathcal L^nh\|_s
		\le C_s\lambda^{-\gamma n}\|h\|_s+C_s|h|_w,
		\qquad h\in\mathcal D,\quad n\ge0.
		\label{strong_stable_estimate}
	\end{equation}
\end{lemma}

\begin{proof}
	Fix $W\in\mathscr W^s$ and
	$\varphi\in C^\gamma(W)$ with
	$|\varphi|_{C^\gamma(W)}\le1$.  For each
	$W_{\mathbf a}\in\mathcal G_n(W)$, choose
	$z_{\mathbf a}\in W_{\mathbf a}$ and set
	\[
	c_{\mathbf a}:=\varphi(F_{\mathbf a}^nz_{\mathbf a}),
	\qquad
	q_{\mathbf a}:=\varphi\circ F_{\mathbf a}^n-c_{\mathbf a}.
	\]
	Stable contraction and the uniform diameter of the stable testing curves
	give
	\[
	|q_{\mathbf a}|_{C^\gamma(W_{\mathbf a})}
	\le C_\gamma \lambda^{-\gamma n}.
	\]
	Combining this with the Lipschitz regularity of the inverse-branch weight
	defined in \eqref{inverse_branch_weight} yields
	\begin{equation*}
		\bigl|q_{\mathbf a}\mathfrak j_{\mathbf a,W}^{(n)}
		\bigr|_{C^\gamma(W_{\mathbf a})}
		\le C_\gamma \lambda^{-\gamma n}
		|\mathfrak j_{\mathbf a,W}^{(n)}|_{\Lip(W_{\mathbf a})},
		\qquad
		|c_{\mathbf a}\mathfrak j_{\mathbf a,W}^{(n)}|_{\Lip(W_{\mathbf a})}
		\le |\mathfrak j_{\mathbf a,W}^{(n)}|_{\Lip(W_{\mathbf a})}.
		\label{stable_test_decomposition_bounds}
	\end{equation*}
	Using the change-of-variables identity \eqref{weak_change_variables}, split
	each pulled-back test function as $q_{\mathbf a}+c_{\mathbf a}$. The first part is controlled by
	$\|h\|_s$ and the constant part by $|h|_w$, and
	Lemma~\ref{lem:inverse_branch_weight_summability} gives
	\[
	\|\mathcal L^nh\|_s
	\le C_\gamma C_J\lambda^{-\gamma n}\|h\|_s+C_J|h|_w.
	\]
	Taking the defining suprema proves \eqref{strong_stable_estimate}.
\end{proof}

\begin{lemma}
	\label{lem:strong_unstable_estimate}
	There is a constant $C_u=C_u(\alpha,\lambda,B_0,\gamma,\beta)$, independent of $n$, such that
	\begin{equation}
		\|\mathcal L^nh\|_u
		\le C_u\lambda^{-\beta n}\|h\|_u
		+C_u\|h\|_s,
		\qquad h\in\mathcal D,\quad n\ge0.
		\label{strong_unstable_estimate}
	\end{equation}
\end{lemma}

\begin{proof}
	Throughout the proof, \(C\) denotes a positive constant that may change from line to line and depends only on the dynamical constants \(\alpha, \lambda, B_0\) and, when interpolation or the anisotropic norms are involved, on the fixed exponents \(\gamma, \beta\). In particular, \(C\) is independent of \(n\), the word \(\mathbf{a}\), the stable leaves \(W_i\), the admissible test functions \(\varphi_i\), and \(h\).
	
	Fix $W_i=W_{t_i}^s$, $i=1,2$, and put
	$
	\delta:=d_{\mathscr W^s}(W_1,W_2)=|t_1-t_2|
	$. 
	Let $\varphi_i\in\Lip(W_i)$ be admissible in the definition of
	$\|\cdot\|_u$.  Thus
	\[
	|\varphi_i|_{\Lip(W_i)}\le1,
	\qquad
	d_{C^\gamma}(\varphi_1,\varphi_2)\le\delta.
	\]
	For a word $\mathbf a$ of length $n$, write
	\[
	W_{\mathbf a,i}:=(F_{\mathbf a}^n)^{-1}
	\bigl(W_i\cap F_{\mathbf a}^n(S_{[\mathbf a]_n})\bigr),
	\qquad 
	\ell_{\mathbf a}:=|I_{\mathbf a}^{(n)}(0)|.
	\]
	Branchwise change of variables gives
	\begin{equation}
		\int_{W_i}\mathcal L^nh\,\varphi_i\,dm
		=\sum_{\mathbf a\in\mathbb{N}^n}
		\int_{W_{\mathbf a,i}}hG_{\mathbf a,i}\,dm,
		\qquad
		G_{\mathbf a,i}:=
		(\varphi_i\circ F_{\mathbf a}^n)
		\mathfrak j_{\mathbf a,W_i}^{(n)}.
		\label{unstable_branch_change_variables}
	\end{equation}
	
	We first record the transverse displacement of one inverse
	branch.  The image under $F_{\mathbf a}^n$ of the bottom interval
	$I_{\mathbf a}^{(n)}\times\{0\}$ is an admissible unstable boundary of
	$F_{\mathbf a}^n(S_{[\mathbf a]_n})$.  The intersections of this boundary with $W_1$ and
	$W_2$ are separated by at most $C\delta$, by (C3) and uniform cone
	transversality.  On the other hand, the derivative of the inverse boundary
	map is $1/J_{\mathbf a,n}^x$, and by \eqref{inverse_horizontal_width_bound}, \(1/J_{\mathbf a,n}^x
	\le C\ell_{\mathbf a}\).
	The inverse crossing points are the intersections of
	$W_{\mathbf a,1}$ and $W_{\mathbf a,2}$ with the bottom transversal.
	Since the leaf metric is exactly the distance between their bottom
	parameters, it follows that
	\begin{equation}
	\delta_{\mathbf a}:=d_{\mathscr W^s}(W_{\mathbf a,1},W_{\mathbf a,2})
	\le C\ell_{\mathbf a}\delta.
	\label{inicial_dis}
	\end{equation}
	
	Next, stable contraction gives
	$|\varphi_i\circ F_{\mathbf a}^n|_{\Lip(W_{\mathbf a,i})}\le1$, while
	\eqref{individual_inverse_weight_bound} gives \(|\mathfrak j_{\mathbf a,W_i}^{(n)}
	|_{\Lip(W_{\mathbf a,i})}
	\le C\ell_{\mathbf a}\). The product estimate consequently yields
	$|G_{\mathbf a,i}|_{\Lip(W_{\mathbf a,i})}\le C\ell_{\mathbf a}$.
	Let $G_{\mathbf a,1}^{\#}$ be the copy of $G_{\mathbf a,1}$ on
	$W_{\mathbf a,2}$, and therefore, their $C^\gamma$ distance is zero.  They are therefore admissible
	for the unstable norm on the pair
	$(W_{\mathbf a,1},W_{\mathbf a,2})$.  Hence
	\begin{equation}
		\left|
		\int_{W_{\mathbf a,1}}hG_{\mathbf a,1}\,dm
		-\int_{W_{\mathbf a,2}}h
		G_{\mathbf a,1}^{\#}\,dm
		\right|
		\le C\ell_{\mathbf a}^{1+\beta}\delta_{\mathbf a}^{\beta}\|h\|_u.
		\label{matched_profile_unstable_bound}
	\end{equation}
	We now compare \(G_{\mathbf a,1}^{\#}\) and \(G_{\mathbf a,2}\).  Write
	\[
	W_{\mathbf a,i}=W_{t_{\mathbf a,i}}^s,
	\qquad
	w_{\mathbf a,i}(y)
	:=\mathfrak j_{\mathbf a,W_i}^{(n)}
	\bigl(\psi(t_{\mathbf a,i},y)\bigr),
	\]
	and define $\Theta_{\mathbf a,i}$ by
	\[
	F_{\mathbf a}^n\bigl(\psi(t_{\mathbf a,i},y)\bigr)
	=\psi\bigl(t_i,\Theta_{\mathbf a,i}(y)\bigr).
	\]
	If $\widetilde\varphi_i=\varphi_i\circ\psi_{t_i}$, then
	\begin{equation}
		\widetilde G_{\mathbf a,i}(y)
		=G_{\mathbf a,i}(\psi(t_{\mathbf a,i},y))
		=\widetilde\varphi_i(\Theta_{\mathbf a,i}(y))
		w_{\mathbf a,i}(y).
		\label{normalized_test_profile_main}
	\end{equation}
	Let us compare \(w_{\mathbf a,1}\) and \(w_{\mathbf a,2}\) through the factors of \eqref{weight_unstable_comparison}. By the same calculation used in Lemma~\ref{lem:transverse_bilip}, the horizontal Jacobians \(J_{\mathbf a,n}^x\)
	differ logarithmically by $O(\delta)$, while (C4) controls the initial and final
	angle factors by $O(\delta^{1/2})$.  Since the cone angles are bounded away from
	zero, it follows that
	\[
	\left|
	\log\frac{w_{\mathbf a,1}(y)}
	{w_{\mathbf a,2}(y)}
	\right|
	\le C\bigl(\delta^{1/2}+\delta\bigr)
	\le C\delta^{1/2}.
	\]
	Finally, by \eqref{individual_inverse_weight_bound},
	\begin{equation}
		\|w_{\mathbf a,1}-w_{\mathbf a,2}\|_\infty
		\le C\ell_{\mathbf a}\delta^{1/2}.
		\label{paired_weight_profile_main}
	\end{equation}
	For the image coordinates, the equal-$y$ points on the two inverse leaves
	are initially separated by $O(\ell_{\mathbf a}\delta)$, by
	\eqref{inicial_dis}.  Their forward images are joined
	in the unstable direction; cylinder distortion and expansion by the factor
	$O(\ell_{\mathbf a}^{-1})$ make the terminal separation $O(\delta)$.
	On the other hand, equal-$y$ points on $W_1$ and $W_2$ are also separated
	by $O(\delta)$ by (C3).  Uniform transversality of the stable and unstable
	cones therefore gives 
	\begin{equation*}
		\|\Theta_{\mathbf a,1}-\Theta_{\mathbf a,2}\|_\infty
		\le C\delta\le C\delta^{1/2}.
		\label{paired_image_coordinate_main}
	\end{equation*}
	Together with
	$|\widetilde\varphi_1-\widetilde\varphi_2|_{C^\gamma}\le\delta$,
	we have
	\[
	\left|
	\widetilde\varphi_1(\Theta_{\mathbf a,1}(y))
	-\widetilde\varphi_2(\Theta_{\mathbf a,2}(y))
	\right|
	\le \delta+C\delta^{1/2}.
	\]
	Expanding the product in \eqref{normalized_test_profile_main}, using
	\eqref{paired_weight_profile_main} and
	$\|w_{\mathbf a,i}\|_\infty\le C\ell_{\mathbf a}$ gives
	\[
	\|\widetilde G_{\mathbf a,1}
	-\widetilde G_{\mathbf a,2}\|_\infty
	\le C\ell_{\mathbf a}\delta^{1/2}.
	\]
	Note that, Lipschitz norms of \(\widetilde G_{\mathbf a,1}\) and \(\widetilde G_{\mathbf a,2}\) are at most \(C\ell_{\mathbf a}\); hence their difference has
	Lipschitz norm at most \(2 C\ell_{\mathbf a}\).  Applying the interpolation inequality
	\[
	|f|_{C^\gamma}
	\le 2\|f\|_\infty^{1-\gamma}\Lip(f)^\gamma
	\]
	and recalling that $\rho=(1-\gamma)/2$, we obtain 
	\begin{equation*}
		\left|
		G_{\mathbf a,1}^{\#}
		- G_{\mathbf a,2}
		\right|_{C^\gamma(W_{\mathbf a,2})}
		\le C\ell_{\mathbf a}\delta^\rho.
		\label{normalized_profile_error}
	\end{equation*}
	Thus 
	\begin{equation}
		\Big|
		\int_{W_{\mathbf a,2}}h
		\bigl( G_{\mathbf a,1}^{\#}
		- G_{\mathbf a,2}\bigr)\,dm
		\Big|
		\le C\ell_{\mathbf a}\delta^\rho\|h\|_s.
		\label{one_leaf_profile_error}
	\end{equation}
	Summing
	\eqref{unstable_branch_change_variables}, and using
	\eqref{matched_profile_unstable_bound} and \eqref{one_leaf_profile_error}, we obtain
	\begin{equation}
		\begin{aligned}
			&\Big|
			\int_{W_1}\mathcal L^nh\,\varphi_1\,dm
			-\int_{W_2}\mathcal L^nh\,\varphi_2\,dm
			\Big|                                                     \\
			&\qquad\le
			C\delta^\beta\|h\|_u
			\sum_{\mathbf a}\ell_{\mathbf a}^{1+\beta}
			+C\delta^\rho\|h\|_s
			\sum_{\mathbf a}\ell_{\mathbf a}.
		\end{aligned}
		\tag{19}
		\label{unstable_sum_before_cylinders}
	\end{equation}
	The intervals $\ell_{\mathbf a}$ have disjoint interiors, so
	$\sum_{\mathbf a}\ell_{\mathbf a}\le1$.  Moreover, the horizontal
	expansion in (H1) and the full-width crossing property give
	$\ell_{\mathbf a}\le \lambda^{-n}$.  Consequently
	\[
	\sum_{\mathbf a}\ell_{\mathbf a}^{1+\beta}
	\le
	\left(\max_{\mathbf a}\ell_{\mathbf a}\right)^\beta
	\sum_{\mathbf a}\ell_{\mathbf a}
	\le \lambda^{-\beta n}.
	\]
	Divide \eqref{unstable_sum_before_cylinders} by $\delta^\beta$ and use
	$\delta^{\rho-\beta}\le1$.  Taking the defining
	suprema proves \eqref{strong_unstable_estimate} on $\mathcal D$.
\end{proof}
The branchwise formula may have finite anisotropic norm without yet being known to belong to the completion of \(\mathcal D\). The following proposition closes this gap: by finite-branch truncation and smoothing, it shows that \(\mathcal L h\) belongs to the corresponding completion for every \(h\in\mathcal D\). The norm estimates then allow \(\mathcal L\) to extend uniquely by density to the completed spaces, see \cite[Lemma~4.1]{De18}, \cite[Proposition~2.7]{DeLi08} and \cite[Lemmas~3.7--3.8]{DeZh11}.

\begin{proposition}
	\label{prop:transfer_well_defined}
	For every $h\in\mathcal D$ and $n\ge1$, the countable branch sum
	$\mathcal L^nh$ belongs to both $\mathcal B_w$ and $\mathcal B$. Moreover $\mathcal L$ extends uniquely to bounded linear
	operators
	\begin{equation}
		\mathcal L:\mathcal B_w\to\mathcal B_w,
		\qquad
		\mathcal L:\mathcal B\to\mathcal B .
		\label{transfer_on_completed_spaces}
	\end{equation}
	The extensions satisfy \eqref{weak_norm_estimate},
	\eqref{strong_stable_estimate} and \eqref{strong_unstable_estimate} for all
	$h$ in the respective space.
\end{proposition}

\begin{proof}
	Fix $h\in\mathcal{D}$ and $n\ge1$.  We first prove that the branchwise sum belongs to the two completed spaces.	We start by the finite-branch truncation. Since the intervals $\ell_{\mathbf a}$ have disjoint interiors, for every
	$\varepsilon>0$ there is a finite set
	$\mathcal A_\varepsilon\subset\mathbb N^n$ such that
	\[
	\sum_{\mathbf a\notin\mathcal A_\varepsilon}
	\ell_{\mathbf a}<\varepsilon.
	\]
	Running the proofs of Lemmas~\ref{lem:weak_norm_estimate},
	\ref{lem:strong_stable_estimate} and \ref{lem:strong_unstable_estimate} over
	the omitted words only, and using
	\eqref{individual_inverse_weight_bound} and
	$\sum_{\mathbf a\notin\mathcal A_\varepsilon}\ell_{\mathbf a}^{1+\beta}
	\le\lambda^{-\beta n}\varepsilon$, gives
	\begin{equation}
	\Bigl\|\sum_{\mathbf a\notin\mathcal A_\varepsilon}
	\mathcal L_{\mathbf a}^nh\Bigr\|_{\mathcal B}
	\le C\varepsilon\,\|h\|_{\mathcal B},
	\qquad
	\Bigl|\sum_{\mathbf a\notin\mathcal A_\varepsilon}
	\mathcal L_{\mathbf a}^nh\Bigr|_w
	\le C\varepsilon\,|h|_w .
	\label{tail_sum}
	\end{equation}
	
	We now smooth each retained branch.	Fix $\mathbf a\in\mathcal A_\varepsilon$.  The density
	$\mathcal L_{\mathbf a}^nh$ is $C^1$ up to the two unstable boundary curves
	of $F_{\mathbf a}^n(S_{[\mathbf a]_n})$.  Extend its smooth factor to a neighborhood of
	the closed strip and replace $\mathbf1_{F_{\mathbf a}^n(S_{[\mathbf a]_n})}$ by a $C^1$
	cutoff whose transition region has width $r$.
	The resulting function $h_{\mathbf a,r}$ belongs to $\mathcal D$.
	
	Let $E_{\mathbf a,r}:=\mathcal L_{\mathbf a}^nh-h_{\mathbf a,r}$.  Uniform
	transversality of stable leaves and unstable boundaries implies that the
	intersection of $\operatorname{supp}E_{\mathbf a,r}$ with any stable leaf has arclength
	$O(r)$.  Hence
	\[
	|E_{\mathbf a,r}|_w+\|E_{\mathbf a,r}\|_s\le C_{\mathbf a,h}r.
	\]
	For two leaves at distance $\delta$, the transverse comparison used in the
	unstable estimate bounds the corresponding difference of cutoff errors by
	$C_{\mathbf a,h}\min\{r,\delta^\rho\}$.  Therefore
	\[
	\|E_{\mathbf a,r}\|_u
	\le C_{\mathbf a,h}
	\sup_{0<\delta\le1}
	\frac{\min\{r,\delta^\rho\}}{\delta^\beta}
	\le C_{\mathbf a,h}r^{1-\beta/\rho}.
	\]
	Because $\beta<\rho$, this tends to zero with $r$.  Thus every retained
	branch density, and hence the finite sum $\sum_{\mathbf a\in\mathcal A_\varepsilon}
	\mathcal L_{\mathbf a}^nh$, belongs to
	both $\mathcal B_w$ and $\mathcal B$.  Letting $\varepsilon\to0$ in
	\eqref{tail_sum} proves
	\[
	\mathcal L^nh\in\mathcal B_w\cap\mathcal B,
	\qquad h\in\mathcal D.
	\]
	
	Finally, The operator is extended to the completions. 
	The inclusion of the Lipschitz unit ball into a fixed multiple of the
	$C^\gamma$ unit ball gives, on $\mathcal D$,
	\begin{equation}
		|h|_w\le C_{ws}\|h\|_s.
		\label{weak_controlled_by_stable}
	\end{equation}
	In particular, the identity on $\mathcal D$ extends continuously from
	$\mathcal B$ to $\mathcal B_w$. Lemma~\ref{lem:weak_norm_estimate} shows that
	$\mathcal L$ maps $|\cdot|_w$-Cauchy sequences in $\mathcal D$ to
	$|\cdot|_w$-Cauchy sequences, which gives the first extension in
	\eqref{transfer_on_completed_spaces}.  Lemmas~\ref{lem:strong_stable_estimate}
	and \ref{lem:strong_unstable_estimate}, together with
	\eqref{weak_controlled_by_stable}, give
	$\|\mathcal Lh\|_{\mathcal B}\le C\|h\|_{\mathcal B}$ on $\mathcal D$ and
	hence the second extension.  Applying the three estimates to seed
	approximations and passing to the limit extends them to the completed
	spaces.
\end{proof}

\subsection{Embedding and compactness}

The next lemma identifies \(\mathcal B_w\) as a space of distributions. The proof relies on the quantitative disintegration established in Lemma~\ref{lem:foliated_Fubini}.

\begin{lemma}
	\label{lem:weak_embedding}
	The identification of $h\in\mathcal D$ with the functional
	\[
	h(\varphi):=\int_Qh\varphi\,dm
	\]
	extends uniquely to a continuous injective linear map
	\begin{equation}
		\mathcal B_w
		\hookrightarrow\bigl(C^1(\mathscr W^s)\bigr)^*.
		\label{weak_embedding}
	\end{equation}
	More precisely,
	\begin{equation}
		\|h\|_{(C^1(\mathscr W^s))^*}
		\le C_{\rm dis}|h|_w,
		\qquad h\in\mathcal B_w.
		\label{weak_embedding_bound}
	\end{equation}
\end{lemma}

\begin{proof}
	Let first $h\in\mathcal D$ and
	$\varphi\in C^1(\mathscr W^s)$.  By
	Lemma~\ref{lem:foliated_Fubini},
	\[
	\begin{aligned}
		\Big|\int_Qh\varphi\,dm\Big|
		=\Big|\int_0^1\int_{W_t^s}
		h\,(\varphi\rho_t)\,dm\,dt\Big|
		\le |h|_w\int_0^1
		|\varphi\rho_t|_{\Lip(W_t^s)}\,dt
		\le C_{\rm dis}|h|_w
		|\varphi|_{C^1(\mathscr W^s)}.
	\end{aligned}
	\]
	Consequently, if $(h_n)\subset\mathcal D$ is Cauchy in the weak norm, then
	$(h_n)$ is Cauchy in
	$(C^1(\mathscr W^s))^*$.  This defines the unique continuous extension and
	proves \eqref{weak_embedding_bound}.
	
	It remains to prove injectivity. Let first $h\in\mathcal D$.  For each
	$t\in[0,1]$, we regard $h$ as a linear functional on
	$\operatorname{Lip}(W_t^s)$ through the identification
	\[
	h(\zeta):=\int_{W_t^s}h\,\zeta\,dm,
	\qquad \zeta\in\operatorname{Lip}(W_t^s).
	\]
	By the definition of the weak norm,
	\(
	|h(\zeta)|
	\le |h|_w\,|\zeta|_{\operatorname{Lip}(W_t^s)},
	\)
	and hence
	\begin{equation}
		|h|_{\operatorname{Lip}(W_t^s)^*}\le |h|_w.
		\label{leafwise_functional_bound}
	\end{equation}
	Thus, for every fixed $t$, this identification extends uniquely and
	continuously from $\mathcal D$ to $\mathcal B_w$.  We use the same symbol
	$h$ for the resulting leafwise functional.  Moreover, the definition of
	$|\cdot|_w$ gives, first on $\mathcal D$ and then by completion,
	\begin{equation}
		|h|_w
		=\sup_{t\in[0,1]}
		|h|_{\operatorname{Lip}(W_t^s)^*},
		\qquad h\in\mathcal B_w.
		\label{weak_norm_as_leafwise_dual_norm}
	\end{equation}
		Suppose now that $h\in\mathcal B_w$ with
	\begin{equation}
		|h|_{C^1(\mathscr{W}^s)^*}=0
		\label{global_functional_zero}
	\end{equation}
	We shall prove that
	\(
	|h|_{\operatorname{Lip}(W_t^s)^*}=0
	\), for all \(t\in [0,1]\). 
	
	Let $\eta\in\operatorname{Lip}([0,1])$.  Transport \(\eta\) horizontally to
	each stable leaf by setting
	\begin{equation*}
		\eta_t(\psi(t,y)):=\eta(y),
		\qquad 0\le y\le1.
		\label{transported_profile}
	\end{equation*}
	Since
	\(
	|y-y'|
	\le d_{W_t^s}\bigl(\psi(t,y),\psi(t,y')\bigr),
	\)
	we have uniformly in $t$,
	\begin{equation}
		|\eta_t|_{\operatorname{Lip}(W_t^s)}
		\le |\eta|_{\operatorname{Lip}([0,1])}.
		\label{transported_profile_bound}
	\end{equation}
	
	Let $A\subset[0,1]$ be an arbitrary Borel set and define
	\begin{equation*}
		\Phi(\psi(t,y))
		:=\mathbf 1_A(t)\,
		\frac{\eta(y)}{\rho_t(\psi(t,y))}.
		\label{localized_global_test}
	\end{equation*}
	The function $\Phi$ is bounded and Borel measurable.  Moreover,
	\eqref{conditional_density_bounds} implies that $1/\rho_t$ is uniformly
	Lipschitz along the leaves; more precisely,
	\[
	\Big|\frac1{\rho_t}\Big|_{\operatorname{Lip}(W_t^s)}
	\le c_{\rm dis}^{-1}
	+C_{\rm dis}c_{\rm dis}^{-2}.
	\]
	Together with \eqref{transported_profile_bound}, this shows that
	\[
	|\Phi|_{C^1(\mathscr W^s)}
	\le C|\eta|_{\operatorname{Lip}([0,1])},
	\]
	where $C$ is independent of $A$.  Hence
	$\Phi\in C^1(\mathscr W^s)$.
	
	Choose $h_n\in\mathcal D$ such that $h_n\to h$ in $\mathcal B_w$.  For
	each $n$, the Fubini and the definition of $\Phi$ give
	\begin{equation}
		\int_Qh_n\Phi\,dm
		=\int_A\int_{W_t^s}h_n\eta_t\,dm\,dt.
		\label{localized_fubini_identity}
	\end{equation}
	For each $n$, the function
	\[
	t\longmapsto h_n(\eta_t)
	\]
	is continuous.  Indeed, after parametrizing $W_t^s$ by
	$y\mapsto\psi(t,y)$, 
	\[
	h_n(\eta_t)=\int_{W_t^s}h_n\eta_t\,dm=\int_0^1 h_n(\psi(t,y))\eta(y)
	\sqrt{1+|\partial_y\chi(t,y)|^2}\,dy,
	\]
	whose dependence on $t$ is continuous by \eqref{global_chart} and
	\eqref{transverse_C4}.  In addition,
	\begin{align*}
		\sup_{t\in[0,1]}
		\left|
		\eta_t(h_n-h_m)
		\right|
		&\le |h_n-h_m|_w
		|\eta|_{\operatorname{Lip}([0,1])}.
	\end{align*}
	Consequently, these continuous functions converge uniformly to
	$t\mapsto h(\eta_t)$.  In particular, this last function is continuous.
		We may now pass to the limit in \eqref{localized_fubini_identity}.  The
	continuity of the global embedding gives
	\[
	\lim_{n\to\infty}\int_Qh_n\Phi\,dm=h(\Phi)=0
	\]
	by \eqref{global_functional_zero}, while the uniform convergence just
	proved gives
	\[
	\lim_{n\to\infty}
	\int_A\int_{W_t^s}h_n\eta_t\,dm\,dt
	=\int_A h(\eta_t)\,dt.
	\]
	Therefore
	\begin{equation*}
		\int_A h(\eta_t)\,dt=0,
				\label{localized_leafwise_vanishing}
	\end{equation*}
    for every Borel set $A\subset[0,1]$. It follows that $h(\eta_t)=0$ for almost every $t$.  Since
	$t\mapsto h(\eta_t)$ is continuous, we conclude that
	\begin{equation}
		h(\eta_t)=0
		\qquad\text{for every }t\in[0,1].
		\label{profile_vanishing_every_leaf}
	\end{equation}
	
	Finally, fix $t_0\in[0,1]$ and let
	$\zeta\in\operatorname{Lip}(W_{t_0}^s)$ be arbitrary.  Define
	\[
	\eta(y):=\zeta(\psi(t_0,y)).
	\]
	Since $\psi_{t_0}$ is Lipschitz, $\eta\in\operatorname{Lip}([0,1])$,
	and by construction $\eta_{t_0}=\zeta$.  Applying
	\eqref{profile_vanishing_every_leaf} to this profile gives
	\[
	h(\zeta)=h(\eta_{t_0})=0.
	\]
	Since $\zeta$ was arbitrary,
	\[
	|h|_{\operatorname{Lip}(W_{t_0}^s)^*}=0.
	\]
	As $t_0$ was arbitrary, \eqref{weak_norm_as_leafwise_dual_norm} yields
	\[
	|h|_w
	=\sup_{t\in[0,1]}
	|h|_{\operatorname{Lip}(W_t^s)^*}
	=0.
	\]
	Thus $h=0$ in $\mathcal B_w$, and the map in
	\eqref{weak_embedding} is injective.
\end{proof}

\begin{lemma}
	\label{lem:compact_embedding}
	The canonical inclusion
	\begin{equation}
		\mathcal B\hookrightarrow\mathcal B_w
		\tag{22}
		\label{strong_weak_compact_embedding}
	\end{equation}
	is compact.
\end{lemma}
\begin{proof}
	Continuity follows from \eqref{weak_controlled_by_stable}.  Since
	$\mathcal B_w$ is complete, it remains to prove that the closed unit ball
	\(\mathcal U:=\{h\in\mathcal B:\|h\|_{\mathcal B}\le1\}\) is totally bounded in $|\cdot|_w$. 
	If $\varphi\in\Lip(W_t^s)$ and
	$|\varphi|_{\Lip(W_t^s)}\le1$, then 
	\(
	\widetilde\varphi:=\varphi\circ\psi_t
	\)
    lies in the fixed set
	\[
	\mathcal K
	:=\{f\in\Lip([0,1]):|f|_{\Lip([0,1])}\le C_s\}.\quad (C_s=\sqrt{1+\alpha^2})
	\]
	Since $0<\gamma<1$, the set $\mathcal K$ is relatively compact in $C^\gamma([0,1])$. 
    Let $\varepsilon>0$.  Choose $r>0$ such that
	\begin{equation}
		2 C_s r^\beta<\frac{\varepsilon}{3},
		\label{compact_leaf_mesh_choice}
	\end{equation}
	and choose a finite $r$-dense $t_1,\ldots,t_M$ of $[0,1]$.  Next choose
	$\kappa>0$ such that
	\begin{equation}
		2 C_s^\gamma\kappa<\frac{\varepsilon}{3},
		\label{compact_test_mesh_choice}
	\end{equation}
	and choose a finite $\kappa$-dense $f_1,\ldots,f_N\in\mathcal K$ of
	$\mathcal K$ in $C^\gamma([0,1])$. For $1\le j\le M$ and $1\le k\le N$, define
	\[
	\Lambda_{j,k}(h)
	:=\int_{W_{t_j}^s}h\,
	(f_k\circ\psi_{t_j}^{-1})\,dm.
	\]
	Each
	$\Lambda_{j,k}$ is a continuous linear functional on $\mathcal B$.
	Consequently, the coordinate image 
	\[
	\bigl\{(\Lambda_{j,k}(h))_{j,k}:h\in\mathcal U\bigr\}
	\subset\mathbb C^{MN}
	\]
	is totally bounded.  Thus, we can write as a finite union $\mathcal{U}:=\cup_{\ell}\{\mathcal{U}_\ell\}$ such that for any $h,h'\in\mathcal{U}_\ell$, 
	\begin{equation*}
		|\Lambda_{j,k}(h-h')|<\frac{\varepsilon}{3}, 
		\label{compact_coordinate_approximation}
	\end{equation*}
	for every $j,k$. Now, let $h,h'\in\mathcal{U}_\ell$, for some $\ell$. Fix a leaf $W_t^s$, and
	$\varphi\in\Lip(W_t^s)$ with $|\varphi|_{\Lip(W_t^s)}\le1$.  Choose
	$t_j$ with $|t-t_j|\le r$, and define $\varphi^\#$ on $W_{t_j}^s$ by 
    \(\varphi^\#\circ\psi_{t_j}=\varphi\circ\psi_t\).
    Hence
	\begin{equation*}
		\Big|
		\int_{W_t^s}(h-h')\varphi\,dm
		-\int_{W_{t_j}^s}(h-h')\varphi^\#\,dm
		\Big|\le
		C_s|t-t_j|^\beta\|h-h'\|_u
		\le2 C_s r^\beta.
		\label{compact_transverse_approximation}
	\end{equation*}
	
	Choose $f_k$ such that
	$|\widetilde\varphi-f_k|_{C^\gamma([0,1])}<\kappa$. 
	It follows from the definition of the strong stable norm that
	\begin{equation*}
		\Big|\int_{W_{t_j}^s}(h-h')
		\bigl(\varphi^\#-f_k\circ\psi_{t_j}^{-1}\bigr)\,dm\Big|
		\le2 C_s^\gamma\kappa.
		\label{compact_test_approximation}
	\end{equation*}
    Adding the last three estimates
	and using \eqref{compact_leaf_mesh_choice}--
	\eqref{compact_test_mesh_choice}, we obtain
	\[
	\Big|\int_{W_t^s}(h-h')\varphi\,dm\Big|<\varepsilon.
	\]
	Taking the supremum over $t$ and over all weak unit tests yields
	$|h-h'|_w<\varepsilon$.
	Choose a finite subset $\mathcal{H}_\varepsilon:=\{h_\ell\in\mathcal{U}_\ell\}$ of $\mathcal{U}$. The conclusion above shows the $\varepsilon$-density of $\mathcal{H}_\varepsilon$ in $\mathcal{U}$ in the $|.|_w$ norm. 
\end{proof}

\subsection{Spectral properties of the transfer operator}
We recall that for isolated point $\lambda\in \sigma(\mathcal{L})$, the \emph{Riesz spectral  
projection} associated with $\lambda$ is defined by 
$$  
\Pi_\lambda:=-\frac{1}{2\pi i}\oint_\Gamma (\mathcal{L}-zI)^{-1}\,dz,
$$  
where $\Gamma$ is any sufficiently small contour enclosing $\lambda$. Subspace $\textrm{Ran}(\Pi_\lambda)$ is called the associated Riesz
spectral subspace. 

Now, Lasota--Yorke inequalities and the compact embedding now give
quasi-compactness.
\begin{theorem}
	\label{thm:quasicompactness}
	The transfer operator $\mathcal L:\mathcal B\to\mathcal B$ is
	quasi-compact.  Its spectral radius is one, and its essential spectral
	radius satisfies
	\begin{equation}
		r_{\rm ess}(\mathcal L|_{\mathcal B})
		\le
		R_*:=\max\{\lambda^{-\gamma},\lambda^{-\beta}\}
		=\lambda^{-\min\{\gamma,\beta\}}<1.
		\label{essential_spectral_radius_bound}
	\end{equation}
	Consequently, for every $\varepsilon>0$, the spectrum in
	$\{z:|z|>R_*+\varepsilon\}$ consists of finitely many isolated eigenvalues,
	each of finite algebraic multiplicity. Equivalently, the associated Riesz
	spectral subspace is finite-dimensional. 
\end{theorem}

\begin{proof}
	We work on the complexification of the spaces.  For $A\ge1$, define the
	equivalent strong norm
	\[
	\|h\|_{\mathcal B,A}:=A\|h\|_s+\|h\|_u.
	\]
	Combining \eqref{strong_stable_estimate} and
	\eqref{strong_unstable_estimate}, for every $n\ge1$ we obtain
	\begin{equation}
		\begin{aligned}
			\|\mathcal L^nh\|_{\mathcal B,A}
			&\le
			\bigl(AC_s\lambda^{-\gamma n}+C_u\bigr)\|h\|_s
			+C_u\lambda^{-\beta n}\|h\|_u
			+AC_s|h|_w\\
			&\le a_{n,A}\|h\|_{\mathcal B,A}+AC_s|h|_w,
		\end{aligned}
		\label{weighted_doeblin_fortet}
	\end{equation}
	where
	\begin{equation*}
		a_{n,A}:=\max\left\{
		C_s\lambda^{-\gamma n}+\frac{C_u}{A},
		C_u\lambda^{-\beta n}
		\right\}.
		\label{weighted_contraction_coefficient}
	\end{equation*}
	The weak estimate gives
	\[
	|\mathcal L^{nk}h|_w\le C_w|h|_w,
	\qquad k\ge0,
	\]
	while Lemma~\ref{lem:compact_embedding} gives the compact embedding of the
	weighted strong space into $\mathcal B_w$.  The standard
	Hennion--Nussbaum argument applied to $\mathcal L^n$ and
	\eqref{weighted_doeblin_fortet} therefore yields
	\begin{equation}
		r_{\rm ess}(\mathcal L^n)\le a_{n,A}.
		\label{hennion_weighted_bound}
	\end{equation}
	
	The essential spectrum is unchanged by replacing the strong norm with an
	equivalent one.  Hence, for fixed $n$, we may let $A\to\infty$ in
	\eqref{hennion_weighted_bound}.  Using
	$r_{\rm ess}(\mathcal L^n)=r_{\rm ess}(\mathcal L)^n$, we get
	\[
	r_{\rm ess}(\mathcal L)
	\le
	\max\left\{
	C_s^{1/n}\lambda^{-\gamma},
	C_u^{1/n}\lambda^{-\beta}
	\right\}.
	\]
	Letting $n\to\infty$ proves
	\eqref{essential_spectral_radius_bound}. So, it remains to identify the spectral radius.  The mass functional
	\[
	\mathfrak m(h):=h(1)
	\]
	is a nonzero continuous functional on $\mathcal B$, and mass preservation
	gives $\mathfrak m\circ\mathcal L=\mathfrak m$.  Thus
	$\mathcal L^*\mathfrak m=\mathfrak m$, so $1\in\sigma(\mathcal L^*)=\overline{\sigma(\mathcal L)}$. Hence, 
	$$r(\mathcal L)\ge1.$$	
	On the other hand, if a spectral value $\zeta$ satisfied $|\zeta|>1$,
	then \eqref{essential_spectral_radius_bound} would make it an eigenvalue of
	finite multiplicity.  For a corresponding nonzero eigenvector $h$,
	the uniform weak estimate \eqref{weak_norm_estimate} would give
	\[
	|h|_w
	=|\zeta|^{-n}|\mathcal L^nh|_w
	\le C_w|\zeta|^{-n}|h|_w,
	\qquad n\ge1.
	\]
	Letting $n\to\infty$ gives $|h|_w=0$. This contradiction give that $r(\mathcal L)\le1$, and hence
	$r(\mathcal L)=1$.  Since $R_*<1$, the operator is quasi-compact.
	
	Finally, the general spectral structure of a quasi-compact operator implies
	that, for every $r>r_{\rm ess}(\mathcal L)$, its spectrum outside
	$\{|z|\le r\}$ is a finite set of isolated eigenvalues of finite algebraic
	multiplicity.  Taking $r=R_*+\varepsilon$ proves the last assertion.
\end{proof}
We complete the spectral picture by identifying the eigenspaces on the unit
circle. We first record the multiplier property needed
to pass between smooth functions and invariant measures with bounded
densities.

\begin{lemma}
	\label{lem:multiplier}
	There is a constant $C_{\rm mult}<\infty$ such that, for every
	$g\in C^1(Q)$ and $h\in\mathcal B$, the product $gh$, initially defined
	on the seed space and then extended by density, belongs to $\mathcal B$ and
	\begin{equation}
		\|gh\|_{\mathcal B}
		\le C_{\rm mult}\,|g|_{C^1(Q)}\|h\|_{\mathcal B}.
		\label{multiplier_bound}
	\end{equation}
	Under the embeddings of Lemma~\ref{lem:weak_embedding}, this extension is
	the usual product of a distribution by a function:
	\[
	(gh)(\varphi)=h(g\varphi).
	\]
	The analogous weak estimate
	\[
	|gh|_w\le C_{\rm mult}|g|_{C^1(Q)}|h|_w
	\]
	also holds for $h\in\mathcal B_w$.
\end{lemma}

\begin{proof}
	We first take $h\in\mathcal D$.  On every $W\in\mathscr W^s$, the
	restriction of $g$ has uniformly controlled Lipschitz norm.  The elementary
	product estimates therefore give
	\[
	|gh|_w\le C|g|_{C^1(Q)}|h|_w, 
	\qquad
	\|gh\|_s\le C|g|_{C^1(Q)}\|h\|_s\quad 
	\]
	It remains to estimate the unstable norm.
	
	Let $W_i=W_{t_i}^s$, put $\delta=|t_1-t_2|\in(0,1]$, and let
	$\varphi_i$ be an admissible pair in the definition of $\|\cdot\|_u$.
	Write
	\(
	\widetilde{g}_i=g\circ\psi_{t_i},
	\widetilde{\varphi}_i=\varphi_i\circ\psi_{t_i},
	\) 
	The $C^0$ transverse Lipschitz estimate for the foliation and the uniform
	leafwise Lipschitz bound give
	\[
	|\widetilde{g}_1-\widetilde{g}_2|_{C^0([0,1])}
	\le C|g|_{C^1(Q)}\delta, 
	\qquad
	\Lip(\widetilde{g}_i)\le C\,|g|_{C^1(Q)}. 
	\]
	Interpolation between the uniform and Lipschitz norms yields
	\begin{equation*}
		|\widetilde{g}_1-\widetilde{g}_2|_{C^\gamma([0,1])}
		\le C|g|_{C^1(Q)}\delta^{1-\gamma}. 
		\label{multiplier_profile_interpolation}
	\end{equation*}
	Moreover, admissibility gives
	$|\widetilde{\varphi}_1-\widetilde{\varphi}_2|_{C^\gamma([0,1])}\le \delta$.  The product estimate in
	$C^\gamma([0,1])$ now implies
	\begin{equation}
		|\widetilde{g}_1\widetilde{\varphi}_1-\widetilde{g}_2\widetilde{\varphi}_2|_{C^\gamma([0,1])}
		\le C|g|_{C^1(Q)}\delta^{1-\gamma}.
		\label{multiplier_product_profiles}
	\end{equation}
	
	Let $\widehat g_2\widehat \varphi_2$ be the horizontal transport of \(g_1\varphi_1\) on $W_2$, that is, \((\widehat g_2\widehat \varphi_2)\circ \psi_2=(g_1 \varphi_1)\circ \psi_1\).
	Both $ g_1 \varphi_1$ and $\widehat g_2\widehat \psi_2$ have Lipschitz norm at most
	$C|g|_{C^1(Q)}$, and they form an admissible pair for the unstable norm. Consequently,
	\[
	\left|
	\int_{W_1}h\,g\varphi_1\,dm
	-\int_{W_2}h\,\widehat g_2\widehat \varphi_2\,dm
	\right|
	\le C|g|_{C^1(Q)}\delta^\beta\|h\|_u.
	\]
	By \eqref{multiplier_product_profiles}, the strong stable norm controls the
	remaining integral:
	\[
	\left|
	\int_{W_2}h\,(\widehat g_2\widehat \varphi_2-g\varphi_2)\,dm
	\right|
	\le C|g|_{C^1(Q)}\delta^{1-\gamma}\|h\|_s.
	\]
	Since $\beta<(1-\gamma)/2<1-\gamma$ and $0<\delta\le1$, the last
	right-hand side is bounded by
	$C|g|_{C^1(Q)}\delta^\beta\|h\|_s$.  Dividing by $\delta^\beta$ and
	taking the defining suprema gives
	\[
	\|gh\|_u
	\le C|g|_{C^1(Q)}(\|h\|_u+\|h\|_s).
	\]
	Together with the stable estimate, this proves
	\eqref{multiplier_bound} on $\mathcal D$.
	
	If $h_n\to h$ in $\mathcal B$, the estimate just proved shows that
	$(gh_n)$ is Cauchy in $\mathcal B$; its limit defines $gh$ and is
	independent of the approximating sequence.  The weak-space statement is
	proved in the same way.  Finally, the identity
	$(gh)(\varphi)=h(g\varphi)$ holds on $\mathcal D$ and passes to the limit
	under the continuous embeddings.
\end{proof}
Let
$
\sigma_{\rm per}(\mathcal L)
:=\{\zeta\in\sigma(\mathcal L):|\zeta|=1\}
$, called the peripheral spectrum of $\mathcal{L}$. For $\zeta\in\sigma_{\rm per}(\mathcal L)$, put $V_\zeta:=\textrm{Ran}(\Pi_\zeta)$. We use
$1$ both for the constant function and for the associated Lebesgue density
on $Q$. As $1$ is an isolated point in $\sigma(\mathcal{L})$, we can define
\begin{equation}
	\mu_0:=\Pi_1 1.
	\label{definition_mu_zero}
\end{equation}
\begin{theorem}
	\label{thm:peripheral_spectrum}
	The peripheral spectrum of $\mathcal L:\mathcal B\to\mathcal B$ has the
	following properties.
	\begin{enumerate}
		\item Every $\zeta\in\sigma_{\rm per}(\mathcal L)$ is semisimple, that is no
		peripheral eigenvalue has a nontrivial Jordan block.
		
		\item The element $\mu_0$ in \eqref{definition_mu_zero} is an
		$F$-invariant Radon probability measure.  Every
		$\nu\in V_\zeta$, for $\zeta\in\sigma_{\rm per}(\mathcal{L})$, is a
		complex Radon measure absolutely continuous with respect to $\mu_0$,
		and
		\[
		\frac{d\nu}{d\mu_0}\in L^\infty(\mu_0).
		\]
		
		\item Every peripheral eigenvalue is a root of unity.  Equivalently, there
		are finitely many integers $q_1,\ldots,q_s\ge1$ such that
		\[
		\sigma_{\rm per}(\mathcal L)
		=\bigcup_{j=1}^s
		\{e^{2\pi i p/q_j}:0\le p<q_j\}.
		\]
		
		\item There are finitely many pairwise mutually singular ergodic
		$F$-invariant probability measures $\mu_1,\ldots,\mu_r$, each
		absolutely continuous with respect to $\mu_0$ with bounded density,
		such that
		\begin{equation*}
			V_1=\operatorname{span}\{\mu_1,\ldots,\mu_r\}.
			\label{invariant_eigenspace_basis}
		\end{equation*}
		That is, these measures form a basis of $V_1$.
	\end{enumerate}
\end{theorem}

\begin{proof}
	The
	Lasota--Yorke estimates, together with $|h|_w\le \|h\|_s$, imply that  
	\begin{equation*}
		\sup_{n\ge0}
		\|\mathcal L^n\|_{\mathcal B\to\mathcal B}
		\le C_{\mathrm{pow}}
		<\infty.
		\label{uniform_strong_power_bound}
	\end{equation*}
	
	We first prove that every peripheral eigenvalue is semisimple. Let $\zeta\in\sigma_{\rm per}(\mathcal L)$, so that $|\zeta|=1$.
	By quasi-compactness, $\zeta$ is an isolated eigenvalue of finite algebraic
	multiplicity.  Suppose, by contradiction, that $\zeta$ is not semisimple.
	Then its generalized eigenspace contains a Jordan chain of length at least
	two.  Without loss the generality, suppose that we have a Jordan chain of length 2, that is, there exist $u,v\in\mathcal B$ such that
	\[
	v\ne0,\qquad
	(\mathcal L-\zeta I)u=v,
	\qquad
	(\mathcal L-\zeta I)v=0.
	\]
	An induction on $n$ gives
	\begin{equation}
		\mathcal L^n u
		=\zeta^n u+n\zeta^{\,n-1}v,
		\qquad n\ge1.
		\label{Jordan_power_formula}
	\end{equation}
	Since $|\zeta|=1$, equation \eqref{Jordan_power_formula} and the uniform
	power bound give
	\[
	n\|v\|_{\mathcal B}
	=\|n\zeta^{\,n-1}v\|_{\mathcal B}
	=\|\mathcal L^n u-\zeta^n u\|_{\mathcal B}
	\le
	\|\mathcal L^n u\|_{\mathcal B}+\|u\|_{\mathcal B}
	\le
	(C_{\mathrm{pow}}+1)\|u\|_{\mathcal B}.
\]
	The right-hand side is independent of $n$, whereas the left-hand side tends
	to infinity because $v\ne0$.  This contradiction shows that every
	$\zeta\in\sigma_{\rm per}(\mathcal L)$ is semisimple. This proves item 1.

 Since the peripheral spectrum is finite and semisimple, for every
	$\zeta\in\sigma_{\rm per}(\mathcal L)$,
	\begin{equation}
		\Pi_\zeta
		=\lim_{N\to\infty}\frac1N
		\sum_{k=0}^{N-1}\zeta^{-k}\mathcal L^k
		\quad\text{in }\mathcal L(\mathcal B,\mathcal B).
		\label{peripheral_cesaro_projector}
	\end{equation}
	The same formula with $\zeta=1$ shows that $\Pi_1$ is a positive operator.  Hence $\mu_0=\Pi_1 1$ is positive, and $\mathcal L\mu_0=\mu_0$ and \(\mu_0(1)=1.\)

Fix $\zeta\in\sigma_{\rm per}(\mathcal L)$ and $\nu\in V_\zeta$.
	The space $V_\zeta$ is finite-dimensional and $\mathcal D$, defined by \eqref{Def:D}, is dense in
	$\mathcal B$, so $\Pi_\zeta\mathcal D=V_\zeta$.  Choose
	$h\in\mathcal D$ with $\nu=\Pi_\zeta h$.  For
	$\varphi\in C^1(Q)$, positivity of the branchwise transfer operator and
	\eqref{peripheral_cesaro_projector} give
	\begin{equation}
		|\nu(\varphi)|
		\le\lim_{N\to\infty}\frac1N
		\sum_{k=0}^{N-1}
		|\mathcal L^kh(\varphi)|                                      
		\le |h|_{C^0(Q)}
		\lim_{N\to\infty}\frac1N
		\sum_{k=0}^{N-1}\mathcal L^k1(|\varphi|)                       
		=|h|_{C^0(Q)}\,\mu_0(|\varphi|).
		\label{Radon_derivative}
	\end{equation}
	In particular,
	$|\nu(\varphi)|\le |h|_{C^0(Q)}|\varphi|_{C^0(Q)}$.  Since \(C^1(Q)\) is uniformly dense in \(C^0(Q)\), $\nu$
	extends uniquely from $C^1(Q)$ to a complex Radon measure on $Q$, and the
	same argument applies to $\mu_0$.  The estimate \eqref{Radon_derivative} shows that \(\nu\) is continuous with respect to the \(L^1(\mu_0)\)-norm. Since \(C^0(Q)\) is dense in \(L^1(\mu_0)\), \(\nu\) extends to a bounded linear functional on \(L^1(\mu_0)\), with norm at most \(|h|_{C^0(Q)}\). By the identification \(\bigl(L^1(\mu_0)\bigr)^*=L^\infty(\mu_0)\),
	there exists \(f_\nu\in L^\infty(\mu_0)\) such that
	\[
	\nu(\varphi)=\int_Q\varphi f_\nu,d\mu_0,
	\qquad
	|f_\nu|_{L^\infty(\mu_0)}
	\le |h|_{C^0(Q)}.
	\]
	Thus \(\nu=f_\nu\mu_0\); in particular, \(\nu\ll\mu_0\) and its Radon--Nikodym density belongs to \(L^\infty(\mu_0)\). This proves item~2.

    We next prove item~3.  We transfer the action of \(\mathcal L\) from measures to their
	densities with respect to \(\mu_0\). Let \(\nu=f\mu_0\in V_\zeta\), where \(f\in L^\infty(\mu_0)\), by item~2.  If \(f\in L^1(\mu_0)\), then
	\(
	\bigl|\mathcal L(f\mu_0)\bigr|
	\le \mathcal L(|f|\mu_0)
	\). 
    We may
	therefore define the relative Perron--Frobenius operator \(P\) by
	\begin{equation}
		\mathcal L(f\mu_0)=(Pf)\mu_0.
		\label{relative_markov_operator}
	\end{equation}
	It is a positive Markov operator on \(L^1(\mu_0)\), with \(P1=1\), and
	satisfies
	\begin{equation}
		\int_Q\varphi\,Pf\,d\mu_0
		=\int_Q(\varphi\circ F)f\,d\mu_0
		\label{markov_koopman_duality}
	\end{equation}
	for every \(f\in L^1(\mu_0)\) and bounded measurable \(\varphi\).
	
	Let \(U\varphi=\varphi\circ F\) be the Koopman operator.  The invariance
	of \(\mu_0\) gives
	\[
	\|U\varphi\|_{L^2(\mu_0)}^2
	=\int_Q|\varphi\circ F|^2\,d\mu_0
	=\int_Q|\varphi|^2\,d(\mathcal L\mu_0)
	=\int_Q|\varphi|^2\,d\mu_0
    =\|\varphi\|_{L^2(\mu_0)}^2
	\]
	so \(U\) is an isometry on \(L^2(\mu_0)\).  Relation
	\eqref{markov_koopman_duality} says that \(P=U^*\) on this space. Now, since
	$
	\mathcal L(f\mu_0)=\zeta f\mu_0,
	$, 
	equation \eqref{relative_markov_operator} gives \(Pf=\zeta f\).  As
	\(|\zeta|=1\), this implies that
	$
	\|Pf\|_{L^2(\mu_0)}=\|f\|_{L^2(\mu_0)}.
	$. 
	Now \(U^*U=I\), while \(UU^*=UP\) is the orthogonal projection onto
	\(\operatorname{Ran}U\).  Hence
	\[
	\begin{aligned}
		\|f\|_{L^2(\mu_0)}^2
		&=\|UPf\|_{L^2(\mu_0)}^2
		+\|(I-UP)f\|_{L^2(\mu_0)}^2\\
		&=\|Pf\|_{L^2(\mu_0)}^2
		+\|(I-UP)f\|_{L^2(\mu_0)}^2.
	\end{aligned}
	\]
	It follows that \(f=UPf\), and therefore
	\begin{equation}
		f=\zeta f\circ F
		\qquad \mu_0\text{-almost everywhere}.
		\label{peripheral_density_cocycle}
	\end{equation}
	Raising \eqref{peripheral_density_cocycle} to the \(k\)-th power and
	using \(PU=U^*U=I\), we obtain
	\begin{equation}
		P(f^k)=\zeta^k f^k,
		\qquad k\ge1.
		\label{powers_of_peripheral_density}
	\end{equation}
	Equivalently,
	\[
	\mathcal L(f^k\mu_0)=\zeta^k f^k\mu_0
	\]
	as finite complex Radon measures.
	
	It remains to verify the nontrivial point that \(f^k\mu_0\) belongs to
	\(\mathcal B\).  The multiplier lemma does not apply directly, since
	\(f^k\) is only known to lie in \(L^\infty(\mu_0)\).  Choose
	\(g_j\in C^1(Q)\) such that
	\[
	\|g_j-f^k\|_{L^1(\mu_0)}\longrightarrow0.
	\]
	By Lemma~\ref{lem:multiplier}, \(g_j\mu_0\in\mathcal B\).  Consider the
	twisted Ces\`aro averages
	\[
	\frac1N\sum_{\ell=0}^{N-1}
	\zeta^{-k\ell}\mathcal L^\ell.
	\]
	Semisimplicity of the peripheral eigenvalues implies that these averages
	converge in operator norm on \(\mathcal B\) to
	\(\Pi_{\zeta^k}\) if \(\zeta^k\) belongs to the peripheral spectrum, and
	to zero otherwise.
	
	On the other hand, \eqref{powers_of_peripheral_density} gives
	\[
	\frac1N\sum_{\ell=0}^{N-1}
	\zeta^{-k\ell}\mathcal L^\ell(f^k\mu_0)
	=f^k\mu_0.
	\]
	Since push-forward by \(F\) is a contraction in total variation,
	\begin{equation}
		\left|
		\frac1N\sum_{\ell=0}^{N-1}
		\zeta^{-k\ell}\mathcal L^\ell(g_j\mu_0)
		-f^k\mu_0
		\right|
		\le
		|(g_j-f^k)\mu_0|
		=\|g_j-f^k\|_{L^1(\mu_0)}.
		\label{projected_smooth_approximation}
	\end{equation}
	If \(\zeta^k\) were not in the peripheral spectrum, letting
	\(N\to\infty\) in \eqref{projected_smooth_approximation}, first against
	\(C^1\) test functions and then by uniform density against all continuous
	test functions, would give
	\[
	\|f^k\mu_0\|_{\rm TV}
	\le \|g_j-f^k\|_{L^1(\mu_0)}.
	\]
	Letting \(j\to\infty\) would imply \(f^k\mu_0=0\), contradicting
	\(f\mu_0\ne0\).  Thus
	\[
	\zeta^k\in\sigma_{\rm per}(\mathcal L).
	\]
	Letting \(N\to\infty\) in
	\eqref{projected_smooth_approximation} now gives
	\[
	\left\|
	\Pi_{\zeta^k}(g_j\mu_0)-f^k\mu_0
	\right\|_{\rm TV}
	\le \|g_j-f^k\|_{L^1(\mu_0)}.
	\]
	Therefore
	\[
	\Pi_{\zeta^k}(g_j\mu_0)
	\longrightarrow f^k\mu_0
	\quad\text{in total variation}.
	\]
	Each term on the left belongs to \(V_{\zeta^k}\).  By item~2,
	\(V_{\zeta^k}\) is a finite-dimensional subspace of the space of finite
	complex Radon measures, and hence it is closed in the total-variation
	norm.  Consequently,
	\[
	f^k\mu_0\in V_{\zeta^k}\subset\mathcal B.
	\]
	
	We have thus proved that \(\zeta^k\) belongs to the peripheral spectrum
	for every \(k\ge1\).  Since the peripheral spectrum is finite, two powers
	of \(\zeta\) must coincide; hence \(\zeta^q=1\) for some \(q\ge1\).
	Thus every peripheral eigenvalue is a root of unity, and the peripheral
	spectrum is a finite union of finite cyclic groups.  This proves item~3.
	
	We finally prove item~4.  Let \(A\subset Q\) be an invariant Borel set,
	that is, \(F^{-1}A=A\), modulo \(\mu_0\). Choose \(g_j\in C^1(Q)\) such that
	\[
	\mu_0(|g_j-\mathbf1_A|)\longrightarrow0.
	\]
	The multiplier lemma gives \(g_j\mu_0\in\mathcal B\).  For
	\(\varphi\in C^1(Q)\), invariance of \(A\) and \(\mu_0\) yields
	\(
	\mu_0(\mathbf1_A\varphi\circ F^\ell)
	=\mu_0(\mathbf1_A\varphi)
	\), for \(\ell\ge0\). It follows that
	\begin{equation*}
		\left|
		(\Pi_1(g_j\mu_0))(\varphi)
		-\mu_0(\mathbf1_A\varphi)
		\right|
		=
		\left|
		\lim_{N\to\infty}\frac1N
		\sum_{\ell=0}^{N-1}
		\mu_0\bigl((g_j-\mathbf1_A)\varphi\circ F^\ell\bigr)
		\right|
		\le
		|\varphi|_{C^0(Q)}\,
		\mu_0(|g_j-\mathbf1_A|).
	\end{equation*}
	By uniform density, this estimate extends to every
	\(\varphi\in C^0(Q)\), and hence
	\[
	\left\|
	\Pi_1(g_j\mu_0)-\mu_0|_A
	\right\|_{\rm TV}
	\le \mu_0(|g_j-\mathbf1_A|)
	\longrightarrow0.
	\]
	Since \(V_1\) is finite-dimensional and therefore closed in total
	variation, we conclude that
	\begin{equation}
		\mu_0|_A\in V_1.
		\label{invariant_restriction_in_eigenspace}
	\end{equation}
	
	Pairwise disjoint invariant sets of positive \(\mu_0\)-measure give,
	through \eqref{invariant_restriction_in_eigenspace}, nonzero mutually
	singular and hence linearly independent elements of \(V_1\).  Since
	\(V_1\) is finite-dimensional, there can be only finitely many such
	sets.  It follows that the invariant sigma-algebra has, modulo
	\(\mu_0\), finitely many atoms \(A_1,\ldots,A_r\).  Define
	\[
	\mu_i:=\frac{\mu_0|_{A_i}}{\mu_0(A_i)},
	\qquad 1\le i\le r.
	\]
	Each \(\mu_i\) is an invariant probability measure belonging to \(V_1\),
	and
	\[
	\frac{d\mu_i}{d\mu_0}
	=\frac{\mathbf1_{A_i}}{\mu_0(A_i)}
	\in L^\infty(\mu_0).
	\]
	Since \(A_i\) is an atom of the invariant sigma-algebra, \(\mu_i\) is
	ergodic.  The measures \(\mu_1,\ldots,\mu_r\) are mutually singular and
	therefore linearly independent.
	
	To prove that they span \(V_1\), let
	\(
	\nu=f\mu_0\in V_1
	\), where 
	\(f\in L^\infty(\mu_0)\) by item~2.  The invariance
	\(\mathcal L\nu=\nu\) gives \(Pf=f\).  The same \(L^2\)-equality argument
	used above, now with \(\zeta=1\), yields
	\[
	f=f\circ F
	\qquad\mu_0\text{-almost everywhere}.
	\]
	An invariant \(L^\infty(\mu_0)\) function is constant almost everywhere
	on each ergodic component \(A_i\).  Thus, for suitable constants
	\(c_1,\ldots,c_r\in\mathbb C\),
	\begin{equation*}
		\nu
		=f\mu_0
		=\sum_{i=1}^r c_i\,\mu_0|_{A_i}
		=\sum_{i=1}^r c_i\mu_0(A_i)\,\mu_i.
	\end{equation*}
	Therefore
	\(
	V_1=\operatorname{span}\{\mu_1,\ldots,\mu_r\}
	\).
	Since these measures are linearly independent, they form a basis of
	\(V_1\).  This proves item~4.
   \end{proof}
\subsection{Simplicity and spectral gap }

Define the quotient
projection along the completed stable leaves by
\begin{equation*}
	q:=\pi_x\circ\psi^{-1}:Q\longrightarrow[0,1],
	\qquad q(W_t^s)=t.
	\label{stable_leaf_quotient_projection}
\end{equation*}

\begin{lemma}
	\label{lem:quotient_measure_ergodic}
	There is a countable full-branch expanding map
	$\bar F:[0,1]\to[0,1]$, defined outside a Lebesgue-null set and satisfying
	\begin{equation*}
		q\circ F=\bar F\circ q
		\qquad \mu_0\text{-almost everywhere},
		\label{quotient_semiconjugacy}
	\end{equation*}
	such that $\bar\mu_0:=q_*\mu_0$ is an invariant probability and
	$(\bar F,\bar\mu_0)$ is exact.  In particular, $\bar F^k$ is ergodic for
	every $k\ge1$.  Moreover,
	\begin{equation}
		d\bar\mu_0=\bar h\,dt,
		\qquad 0<c_\mu\le \bar h\le C_\mu<\infty
		\quad\text{almost everywhere}.
		\label{quotient_density_bounds}
	\end{equation}
\end{lemma}

\begin{proof}
	Each post $S_i$ is saturated by the completed stable foliation, so that
	$\psi^{-1}(S_i)=I_i\times[0,1]$ for an interval $I_i$.  Branchwise
	invariance of stable leaves defines
	$\bar F_i:I_i\to[0,1]$ by
	$F_i(W_t^s)\subset W_{\bar F_i(t)}^s$. The full-width property of $U_i$
	makes every $\bar F_i$ onto. The maps \(\bar F_i\) therefore define a
	countable full-branch map \(\bar F\) on \([0,1]\) outside the endpoints 
	of the intervals \(I_i\), and by construction \(q\circ F=\bar F\circ q\).
	
	For \(n\ge1\) and \(\mathbf a\in\mathbb N^n\), recall that
	\(I_{\mathbf a}^{(n)}(0)=\{x:(x,0)\in S_{[\mathbf a]_n}\}\). 
	The restriction	of \(\bar F^n\) to \(I_{\mathbf a}^{(n)}(0)\) 
	is a full branch of \(\bar F^n\). We denote its inverse by
	\(g_{\mathbf a}^{(n)}\). The inverse-transverse weight estimates
	\eqref{individual_inverse_weight_bound} and
	\eqref{weight_unstable_comparison}, together with the
	horizontal distortion estimate
	\eqref{horizontal_cylinder_distortion} and the two-sided density
	bounds of Lemma~\ref{lem:foliated_Fubini}, imply that there exists
	\(D_q\ge1\), independent of \(n\) and \(\mathbf a\), such that
	\begin{equation}
		D_q^{-1}|I_{\mathbf a}^{(n)}(0)|
		\le
		|(g_{\mathbf a}^{(n)})'(t)|
		\le
		D_q|I_{\mathbf a}^{(n)}(0)|
		\quad\text{for a.e. }t\in[0,1].
		\label{quotient_inverse_distortion_main}
	\end{equation}
	In particular, the \(n\)-step branch intervals
	\(I_{\mathbf a}^{(n)}(0)\) shrink uniformly to points as
	\(n\to\infty\).
	
	Let \(\bar m:=q_*m\). By the foliated Fubini formula, \(d\bar m=r_0(t)\,dt\), 
	where \(r_0\) is bounded above and below by positive constants.
	Since the intervals \(I_{\mathbf a}^{(n)}(0)\),
	\(\mathbf a\in\mathbb N^n\), partition \([0,1]\) modulo endpoints,
	and their lengths sum to one estimate \eqref{quotient_inverse_distortion_main} 
	gives uniform two-sided bounds for the densities of
	\((\bar F^n)_*\bar m\), independently of \(n\).
	
	The Ces\`aro averages of \(\mathcal L^n1\) converge to \(\mu_0\).
	Pushing these averages forward by \(q\) and using
	\(q\circ F=\bar F\circ q\), therefore proves \eqref{quotient_density_bounds} and
	\(\bar F\)-invariance of $\bar\mu_0$. 
	
	The distortion estimate
	\eqref{quotient_inverse_distortion_main} and the two-sided bounds
	on \(\bar h\) imply that there exists \(K_q\ge1\) such that, for
	every \(n\ge1\), every \(\mathbf a\in\mathbb N^n\), and every
	measurable \(B\subset[0,1]\),
	\begin{equation}
		K_q^{-1}\bar\mu_0(B)
		\le
		\frac{
			\bar\mu_0(g_{\mathbf a}^{(n)}(B))
		}{
			\bar\mu_0(I_{\mathbf a}^{(n)}(0))
		}
		\le
		K_q\bar\mu_0(B).
		\label{quotient_relative_measure_distortion_main}
	\end{equation}
	Let $A$ belong to the tail sigma-algebra
	$\bigcap_{n\ge0}\bar F^{-n}\mathscr B([0,1])$ and suppose
	$\bar\mu_0(A)>0$. For every \(n\ge1\), choose a
	measurable set \(A_n\) such that \(A=\bar F^{-n}A_n
	\) modulo \(\bar\mu_0\). Since \(\bar\mu_0\) is equivalent to 
	Lebesgue measure, we may choose	a Lebesgue density point \(t\in A\).
	
	For each \(n\), let \(I_{\mathbf a}^{(n)}(0)\) be the \(n\)-step
	branch interval containing \(t\), where \(\mathbf a\) depends on
	\(n\). Since these intervals shrink to \(t\), the Lebesgue
	differentiation theorem and the two-sided bounds on \(\bar h\) give
	\[
	\frac{
		\bar\mu_0(I_{\mathbf a}^{(n)}(0)\setminus A)
	}{
		\bar\mu_0(I_{\mathbf a}^{(n)}(0))
	}
	\longrightarrow0.
	\]
	On \(I_{\mathbf a}^{(n)}(0)\), the relation
	\(A=\bar F^{-n}A_n\) gives
	\[
	I_{\mathbf a}^{(n)}(0)\setminus A
	=
	g_{\mathbf a}^{(n)}([0,1]\setminus A_n)
	\qquad\text{modulo null sets}.
	\]
	Applying \eqref{quotient_relative_measure_distortion_main} therefore
	yields
	\[
	\bar\mu_0([0,1]\setminus A_n)
	\le
	K_q
	\frac{
		\bar\mu_0(I_{\mathbf a}^{(n)}(0)\setminus A)
	}{
		\bar\mu_0(I_{\mathbf a}^{(n)}(0))
	}
	\longrightarrow0.
	\]
	On the other hand, invariance of \(\bar\mu_0\) and
	\(A=\bar F^{-n}A_n\) imply \(\bar\mu_0(A)=\bar\mu_0(A_n)\) 
	for every \(n\). Consequently \(\bar\mu_0(A)=1\), proving that
	\((\bar F,\bar\mu_0)\) is exact.
	
	Finally, since the sigma-algebras
	\(\bar F^{-n}\mathscr B([0,1])\) form a decreasing sequence, for
	every \(k\ge1\),
	\[
	\bigcap_{j\ge0}
	(\bar F^k)^{-j}\mathscr B([0,1])
	=
	\bigcap_{j\ge0}
	\bar F^{-kj}\mathscr B([0,1])
	=
	\bigcap_{n\ge0}
	\bar F^{-n}\mathscr B([0,1]).
	\]
	Thus every power \(\bar F^k\) is exact and, in particular, ergodic.
\end{proof}

\begin{lemma}
	\label{lem:lift_ergodicity_simplicity}
	The system $(F,\mu_0)$ is totally ergodic: $(F^k,\mu_0)$ is ergodic for
	every $k\ge1$.  Consequently,
	\begin{equation}
		V_1=\ker(\mathcal L-I)=\mathbb C\mu_0,
		\label{simple_invariant_eigenspace}
	\end{equation}
	and the eigenvalue $1$ of $\mathcal L:\mathcal B\to\mathcal B$ is
	simple.
\end{lemma}

\begin{proof}
	Fix $k\ge1$ and let $A$ be an invariant set for $F^k$ with $\mu_0(A)>0$. Then $\bar{A}:=q(A)$ is an invariant set for $\bar{F}^k$. 
    Since $\bar{\mu_0}$ is an  ergodic measure for $\bar{F}^k$, 
    $\bar\mu_0(\bar{A})=1$. Thus
    $$\mu_0(A)=\lim_{n\to\infty} \mu_0\big(F^{kn}(q^{-1}(\bar{A})\big)=\bar\mu_0(\bar{A})=1.$$
	
	If $\nu\in V_1$, Theorem~\ref{thm:peripheral_spectrum} gives
	$\nu=f\mu_0$ with $f\in L^\infty(\mu_0)$.  The Markov-operator argument
	in its proof gives $f\circ F=f$.  Ergodicity makes $f$ constant, proving
	\eqref{simple_invariant_eigenspace}.  Item~1 of
	Theorem~\ref{thm:peripheral_spectrum} says that $1$ is semisimple, so its
	geometric and algebraic multiplicities are both one.
\end{proof}

\begin{theorem}
	\label{thm:spectral_gap}
	The only spectral value of $\mathcal L:\mathcal B\to\mathcal B$ on the
	unit circle is the simple eigenvalue $1$.  Consequently, there exist a
	bounded operator $\mathcal N:\mathcal B\to\mathcal B$, constants
	$C_{\rm gap}<\infty$ and $\vartheta\in(0,1)$, and the rank-one projector
	\begin{equation}
		\Pi_1h=h(1)\mu_0
		\label{spectral_gap_rank_one_projector}
	\end{equation}
	such that
	\begin{equation}
		\mathcal L^nh=\Pi_1h+\mathcal N^nh,
		\qquad
		\|\mathcal N^n h\|_{\mathcal B}
		\le C_{\rm gap}\vartheta^n\|h\|_{\mathcal B},
		\qquad n\ge1.
		\label{spectral_gap_decomposition}
	\end{equation}
	In particular,
	\[
	\sigma(\mathcal L)
	\subset\{1\}\cup\{z\in\mathbb C:|z|\le\vartheta\}.
	\]
\end{theorem}

\begin{proof}
	Let $\zeta\in\sigma_{\rm per}(\mathcal L)$ and choose
	$0\ne\nu=f\mu_0\in V_\zeta$ as in
	Theorem~\ref{thm:peripheral_spectrum}.  Equation
	\eqref{peripheral_density_cocycle} gives
	$f=\zeta f\circ F$, and item~3 of that theorem says that $\zeta$ is a
	root of unity.  If $\zeta^k=1$, iteration gives $f=f\circ F^k$.
	Total ergodicity makes $f$ constant; since $f\ne0$, the original identity
	then forces $\zeta=1$.  Hence
	$\sigma_{\rm per}(\mathcal L)=\{1\}$, and the preceding lemma shows that
	$1$ is simple.
	
	Quasi-compactness and $r_{\rm ess}(\mathcal L)<1$ now imply that the
	spectrum different from $1$ is contained in a disk of radius strictly
	smaller than one.  Let $\Pi_1$ be the Riesz projector at $1$ and put
	$\mathcal N=\mathcal L-\Pi_1$.  Since
	$V_1=\operatorname{Ran}\Pi_1=\mathbb C\mu_0$, the Ces\`aro formula for
	$\Pi_1$, mass preservation, and $\mu_0(1)=1$ give
	\eqref{spectral_gap_rank_one_projector}.  Moreover,
	$\Pi_1\mathcal N=\mathcal N\Pi_1=0$ and $r(\mathcal N)<1$.  Choosing
	$\vartheta\in(r(\mathcal N),1)$, the spectral-radius formula gives
	$\|\mathcal N^n\|\le C_{\rm gap}\vartheta^n$.  Since
	$\mathcal L=\Pi_1+\mathcal N$, this proves
	\eqref{spectral_gap_decomposition}.
\end{proof}

\subsection{Limit theorems}

The spectral gap established above implies exponential decay of
correlations. Moreover, the multiplier estimate of
Lemma~\ref{lem:multiplier} allows the perturbative argument for twisted
transfer operators to be carried out on \(\mathcal B\). We obtain the
following standard statistical consequences. The proof follows the
argument of Section~4.3 of \cite{De18}, based on analytic perturbation
theory for isolated eigenvalues; see also \cite{Kato}.

\begin{corollary}
	\label{thm:statistical_limit_laws}
	The invariant probability measure \(\mu_0\) satisfies the following
	properties.
	\begin{itemize}
		\item \emph{Exponential decay of correlations.}
		There exist constants \(C_{\rm cor}<\infty\) and
		\(\vartheta_{\rm cor}\in(0,1)\) such that, for every
		\(f,g\in C^1(Q)\) and every \(n\ge0\),
		\begin{equation*}
			\left|
			\int_Q f\,g\circ F^n\,d\mu_0
			-\int_Q f\,d\mu_0\int_Q g\,d\mu_0
			\right|
			\le
			C_{\rm cor}\vartheta_{\rm cor}^n
			|f|_{C^1(Q)}|g|_{C^1(Q)} .
			\label{exponential_correlation_decay}
		\end{equation*}
		
		\item \emph{Central limit theorem.}
		Let \(\varphi\in C^1(Q)\) be real-valued and satisfy
		\(\mu_0(\varphi)=0\). Then there exists
		\(\sigma_\varphi^2\ge0\) such that
		\begin{equation*}
			\frac{1}{\sqrt n}
			\sum_{j=0}^{n-1}\varphi\circ F^j
			\xrightarrow{\;\text{dist.}\;}
			\mathcal N(0,\sigma_\varphi^2).
			\label{central_limit_theorem}
		\end{equation*}
	\end{itemize}
\end{corollary}
\section{Virtually Expanding and Existence of ACIP}\label{SecV.E.}
Here, we suppose the infinite smoothness of the branch maps $F_i$. Let us outline the proof, which closely follows the strategy of Tsujii \cite{T23}. In \cite{T23}, the author proves the Lasota--Yorke inequality in the Sobolev space $H^\mu$ with $\mu > 0$ first by proving it for smooth functions supported in small open sets and then extending the result globally via a partition of unity. For our local version, we adopt the same initial setup; however, for globalization, we require the Sobolev space to be close under the multiplication by characteristic functions. This holds in $H^\mu$ only for $\mu < 1/2$, thus determining our choice of space.

The questions that can be raised next is whether a method can be devised in which a higher regularity is achieved? Since a GHM has a Markov structure, meaning that all discontinuity points are mapped to boundary points, we believe that this question is not far-fetched.
\begin{figure}[H]
	\begin{center}
\includegraphics[height=5cm,width=11cm]{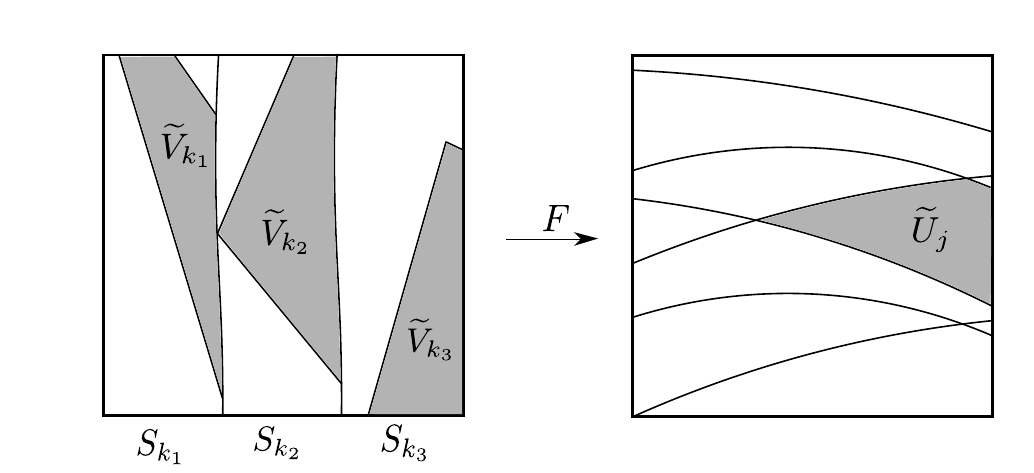}
\caption{Local Mechanism for Lasota York Inequality\label{VirtuallyExpanding}  }
\end{center}
\end{figure}
\subsection{Local Lasota-York Inequality}
For the local version, we need some notation. Suppose that $N$ strips and put $U:=\cup_{i=1}^N U_i$. Then, $U$ can be partitioned as 
$$U:=\bigcup_{j=1}^{\mathrm{N}(F)} \widetilde{U}_j$$ 
such that for each $j$, $\widetilde{U}_j$ is an open subset of $U_i$, for some $i$, and there is a finite set $\Uptheta_j$ of natural numbers with $\#\,\Uptheta_j<\mathrm{N}(F)$ such that for any $k\in \Uptheta_j$, there are open sets $\widetilde{V}_k\subset S_k$ with $F_k:\widetilde{V}_k\to \widetilde{U}_j$ is a $\mathcal{C}^\infty$ smooth diffeomorphism  (see Figure \ref{VirtuallyExpanding}). For any $j$, put 
$$\widetilde{W}_j:=\bigcup_{k\in\Uptheta_j} \widetilde{V}_k.$$
We first state the local Lasota-York inequality as described in the Apprendix II.
\begin{theorem}[Appendix II]\label{T24}
Let $F$ be $\mu$-virtually expanding, $\max\{0,\mu-1/2\}<\mu^{'}<\mu$ and $u$ a $C^\infty$ map supported on $\widetilde{W}_j$. Then, there is a constant $C>0$ such that
$$\|\mathcal{L}u\|\leq 
\sqrt{B^{2\mu}(F)+\epsilon}\|u\|_{H^\mu}+
C\|u\|_{H^{\mu'}}
$$
\end{theorem}
\subsection{Global Lasota-Yorke inequality}
We now globalize the local estimate. For this, we need the following results saying that for $\mu<1/2$, Sobolev space $H^\mu$ is closed with respect to 
the multiplication of characteristic functions.
\begin{theorem}\label{Regularity of Sobolev}
Let $0 \le \mu < \frac12$ and let $E \subset \mathbb{R}^2$ be a set with piecewise $C^1$ boundary. Define the multiplication operator
$
\mathcal{M}: H^\mu(\mathbb{R}^2) \to H^\mu(\mathbb{R}^2)$, defined by 
$\mathcal{M}u = \chi_E \, u.
$
Then 
\[\lim_{\mu\to 0}\|\mathcal{M}\|_{H^\mu} =1.
\]
\end{theorem}
\begin{proof}
It is known that \(H^\mu(\mathbb{R}^2)=[L^2(\mathbb{R}^2),H^{\mu_0}(\mathbb{R}^2)]_\theta\) with \(\theta:=\mu/\mu_0\) for any fixed \(\mu_0\in(\mu,1/2)\) (see \cite[Chapter 4]{Tay11}). Now, since \(\|\mathcal{M}\|_{L^2}=1\), the interpolation theorem for operator norms yields
\[
1\leq \|\mathcal{M}\|_{H^\mu}\le\|\mathcal{M}\|_{H^{\mu_0}}^{\theta}.
\]
This proves the claimed limit.
\end{proof}
Now, suppose that a smooth function $u$ defined on $S$. Put $u=\sum_{j=1}^{\mathbf{N}(F)} \chi_ju$, where $\chi_j:=\chi_{\widetilde{W}_j}$ is the characteristic function of $\widetilde{W}_k$. 
By iterating the inequality above, we get
$$\|\mathcal{L}^n(u)\|_{H^\mu}\leq \|\mathcal{M}\|_{H^\mu}\big(\sqrt{B^{2\mu}(F^n)+\epsilon}\,\,\mathrm{N}(F^n)\|u\|_{H^\mu}+C\,\mathrm{N}(F^n)\|u\|_{H^{\mu^\prime}}
\big).$$
Now, using Theorem \ref{Regularity of Sobolev}, and taking $n$ large enough so that $\sqrt{B^{2\mu}(F^n)+\epsilon}\,\,\mathrm{N}(F^n)<\bar{\lambda}<1$, we get the desired Lasota-York inequality for some $\mu<1/2$. That is
$$\|\mathcal{L}^n(u)\|_{H^\mu}\leq \bar{\lambda}\|u\|_{H^\mu}+C^*\|u\|_{H^{\mu^\prime}}
.$$

\begin{example}\label{Example}
Let $n\geq 2$ be a natural number and choose $\lambda$ with  
$1/n<\lambda\leq \sqrt{1/n}<1$. Put $\theta_0:=\arctan(1-\lambda)$.  For $j=0,\ldots, n-1$, let $S_j:=[j/n,(j+1)/n]\times [0,1]$ and 
$\theta_j:=\theta_0\big(1-\frac{2j}{n-1}\big)$. Let $F_j$ maps $S_j$ linearly to a parallelogram of sides $n$ and $\lambda$ and angel $\theta_j$ (see Figure \ref{Example}). That is 
$$F_j(x,y)=\big(nx-j,\lambda y+\tan\theta_j(nx-j)+\frac{j(1-\lambda)}{n-1}\big).$$
\begin{figure}[H]
	\begin{center}
\hspace{1.3cm}\includegraphics[height=5cm,width=10.5cm]{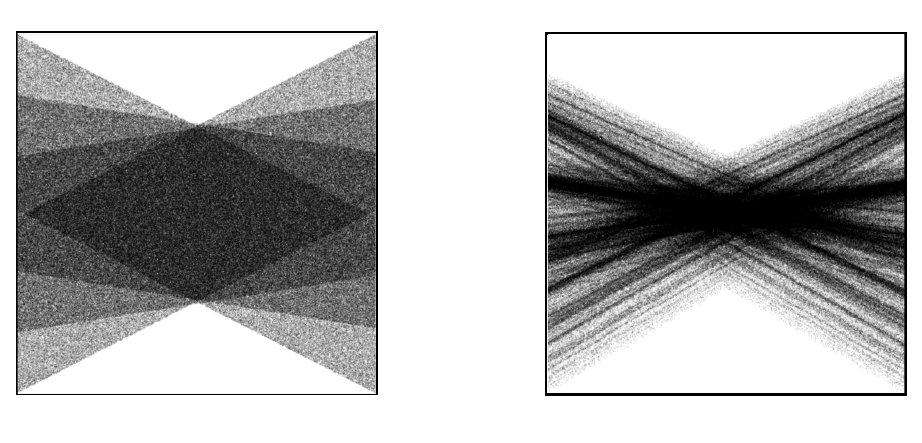}
\caption{Mapping $F$ in Example \ref{Example} for $n=4$, $\lambda=0.49$. The left frame is 
the first iteration and the right frame is tenth iteration. $8.10^5$ points are shown. The asymmetry in the iterations $>1$ arises from the fact that all the parallelograms are from the left boundary to the right boundary. }
\end{center}
\end{figure}
It is clear that
$\mathrm{Jac}(F)=\lambda n$. A more subtle observations shows that for any natural number $l$, 
$\mathrm{N}(F^l)\leq \lambda^{l-1} n^l$.
Now, put  $
N:=[\,\lambda^{l-1}n^{l}\,]+1
$ and 
$a_j:=\theta_0(1-\frac{2j}{N})$, for $j=0,1,\dots,N$. A direct calculating gets that for any $(x,v)\in T^1S$ and iteration $F^l$,
\begin{align*}
b^{1/2}(x,v)&\leq\frac{1}{(n\lambda)^{l}}
\sum_{j=0}^{N}
\frac{1}{\sqrt{n^{2l}\sin^2 a_{j}+\lambda^{2l}\cos^2 a_{j}}}\\
&\simeq\frac{1}{(n\lambda)^{l}}\Big(\frac{N}{2\theta_0}\int_{-\theta_0}^{\theta_0}
\frac{\mathrm{d}\theta}{\sqrt{n^{2l}\sin^2 \theta+\lambda^{2l}\cos^2 \theta}}\Big)\\
&\le \frac{1}{(n\lambda)^{l}} \Big(\frac{\lambda^{l-1}}{\theta_0}\log(\frac{n}{\lambda})+\mathcal{O}(\lambda^{l-l})\Big)
\\&\leq \frac{C}{n^l}\log(\frac{n}{\lambda}).
\end{align*}
In particular, 
$$\sqrt{B^1(F^l)}\,\mathrm{N}(F^l)\leq C (\lambda\sqrt{n})^l.$$
That means in particular that $F$ admits a unique absolutely continuous invariant measure. 
\end{example}
\section{Appendix I (Virtual Expanding)}\label{Ap1}
 Let $M$ be a closed 
smooth Riemannian manifold.
Let also $F:M\to M$ be a $C^{1}$ endomorphism and the Jacobian of $f$ doesn't vanish at any point of the manifold. We denote by $T^*M$ the cotangent bundle of $M$ and define \(F^{\dag}:T^*M\to T^*M\) by
\[F^{\dag}\big((q, \zeta); F\big):=\big(F(q), \big((D_{q}F)^{*}\big)^{-1}(\zeta)\big),\]
where $\zeta$ is a linear functional on $T_q M$ and consequently,  $((D_{q}F)^{*})^{-1}(\zeta)$ is a linear functional on $T_{F(q)}M$. For every nonzero linear functional  $(p, \eta)$ on $T^*M$ and every positive real number $\mu$, put 
\[b^{\mu}((p, \eta) ; F):=\sum_{F^{\dag}(q,\zeta )=(p, \eta)} \frac{(\|\eta\|/\|\zeta\|)^{\mu}} {|J F(q)|},\] 	
where $\|.\|$ is the operator norm associated to the Riemannian metric on $T^*M$.
Now, we give an equivalent geometric definition of virtually expanding endomorphisms on the tangent bundle. 
which is derived from the Riesz representation theorem. First, assume that $p\in M$ has $n$ pre-images $ q_1,q_2,...,q_n$;  and $\zeta_i:=(D_{q_i}f)^* \eta$  which is defined by $\zeta_i(w):=(D_{q_i}f)^* \eta(w)=\eta(D_{q_i}f(w))$. For $(p,v)\in T^{1}M$, let $\pi_{v}:T_{p}M\to \langle v\rangle$ denote the orthogonal projection onto the one-dimensional subspace generated
by $v$.
For every pre image $q_i$ of $p$, define 
$\eta(q_i,v):=\pi_v\big((D_{q_i}F(B_1(q_i))\big)$, where 
$B_1(q_i):=\{w\in T_{q_i}M:\|w\|\le1\}$ is the unit ball of $T_{q_i}M$. Observe that
\[\eta(q_i,v)=\sup_{\|w\|\le1}\left\|\pi_v(D_{q_i}F(w))\right\|.\]
Thus, by using the Riesz representation theorem. then we have,
\begin{align*}
	 b^{\mu}((p, \eta) ; F)& =\sum_{F^{\dag}(q, \zeta)=(p, \eta)}\frac{(\|\eta\| /\|\zeta\|)^{\mu} }{|J F(q) |}=\sum_{i=1}^{n}\frac{\big(\|\eta\|/\| \zeta_{i} \|\big)^{\mu}}{\big|J F\left(q_{i}\right)\big|}
	\\
 & =\sum_{i=1}^{n}\frac{1}{\left|J F\left(q_{i}\right)\right|} \Big(\sup _{\substack{|v|\le1 \\
			v \in T_{p} M}}|\eta(v)|\Big)^{\mu} \Big( \sup _{\substack{|w|\le1 \\
			w \in T_{q_{i}} M}}\left|\eta\left(D_{q_{i}} F(w)\right)\right|\Big)^{-\mu} \\
 & =
\sum_{i=1}^{n}
\frac{1}
{|JF(q_i)|\,
\eta(q_i,v)^{\mu}}.
 \end{align*}

\section{Appendix II (Local Lasota-York Inequality for Virtually Expanding GHMs)}

\subsection{Pseudo differential operators and Sobolev spaces}
Let $M$ be a close manifold and \(U\subset M\) an open set. For \(m\in\R\), the symbol class \(S^m_{1,0}(U)\) consists of all \(C^\infty\) functions \(a(x,\xi)\) on \(U\times\R^2\) such that, for every pair of multi indices \(\alpha,\beta\), and every compact $K\subset U$, there is a constant \(C_{K,\alpha,\beta}\) with
\[
   |\partial_x^\alpha\partial_\xi^\beta a(x,\xi)|
   \leq C_{K,\alpha,\beta}\ang{\xi}^{m-|\beta|}
   \qquad (x\in K,\ \xi\in\R^2),
\]
where
\[
   \ang{\xi}:=(1+|\xi|^2)^{1/2}.
\]
For any $a\in S^m_{1,0}(U)$, we can associate the pseudo differential operator as follow,
\begin{equation}\label{PsudoO}
   \big(\operatorname{Op}(a)u\big)(x)
   :=(2\pi)^{-d}\int_{\R^2}\int_{\R^2}
       e^{i(x-y)\cdot\xi}a(x,\xi)u(y)\,dy\,d\xi,
\end{equation}
where  \(u\in C_c^\infty(U)\).  We write \(\Psi^m(M)\) for the class of pseudodifferential operators of order \(m\). For \(s\in\R\), \(H^s(M)\) is the usual Sobolev space of order \(s\).  For \(s\geq0\) we use the globally defined elliptic weight
\begin{equation*}\label{DefW}
   W^s(x,\xi):=\ang{|\xi|_x}^{s}=(1+|\xi|_x^2)^{s/2},
\end{equation*}  
and use the equivalent norm
\[
   \|u\|_{H^s}^2
   :=\|\op(W^s)u\|_{L^2}^2+\|u\|_{L^2}^2.
\]
When a negative order appears below, we use the convention that \(H^s(M)\) is the usual \(L^2(M)\) space. It is known that for $\mu'\geq \mu$,  
\begin{equation*}
\label{Sineq}\|u\|_{H^{\mu'}}\leq C(\mu,\mu')
\|u\|_{H^\mu}
\end{equation*} 
\begin{lemma}[Basic properties on pseudodifferential operators]\label{BPDO}
The following facts will be used.
\begin{enumerate}
\item \textbf{Boundedness theorem.} If \(A\in\Psi^m(M)\), then, for every \(s\in\R\), $A:H^s(M)\to H^{s-m}(M)$ continuously.
\item If \(A=\op(a)\in\Psi^m\) and \(B=\op(b)\in\Psi^n\), then
$AB=\op(ab)+R$ where $R\in\Psi^{m+n-1}$. Moreover, $AB$ is again a pseudo differential operator, that is 
$AB = \operatorname{Op}(c)$, for some symbol $c(x,\xi)$.
\item If \(A=\operatorname{Op}(a)\in\Psi^m\), then
$ A^*=\operatorname{Op}(\overline a)+R$, for some $R\in\Psi^{m-1}$. In particular, 
if \(a\) is real-valued, then \(A^*-A\in\Psi^{m-1}\).
\item If \(A\in\Psi^m\) and \(B\in\Psi^n\), then the commutator satisfies
$[A,B]=AB-BA\in\Psi^{m+n-1}$, 
provided the principal scalar symbols commute, which they do for scalar symbols.
\end{enumerate}
\end{lemma}
\begin{proof} (see \cite[Chapter II, \S4]{Tay81}) .
\end{proof}
\begin{lemma}[Sharp G\aa rding inequality in the form used here]
Let $p(x,\xi) \in S_{1,0}^m(U)$, where $0 \le \delta < \rho \le 1$. If $\operatorname{Re} p(x,\xi) \ge 0$ 
then there exists a constant $C$ such that for all $u \in C_0^\infty(U)$,
\begin{equation*}
\operatorname{Re} (\op(p)u, u) \ge - C \|u\|_{H^{(m-1)/2}}^2,
\end{equation*}
where $\op(p)$ is the pseudodifferential operator associated to $p$ by (\ref{PsudoO}).
\end{lemma}
\begin{proof}  (see \cite[Theorem 2.3]{Tay81})
\end{proof}
\begin{lemma}[Compact Sobolev embedding]
If \(s>t\), then the inclusion
$ H^s(M)\hookrightarrow H^t(M)$ is compact.
\end{lemma}
\begin{proof}
    (see \cite[Chapter 4, Section 4.3, Proposition 3.4.]{Tay11})
\end{proof}




\subsection{Local Lasota-Yorke inequality}
Here, we essentially follow the proof of the Lasota–York inequality in a local form, as given in \cite{T23}. We use the notation of Section \ref{SecV.E.}. Let
$\varphi_k=F_k^{-1}:\widetilde{U}_j\to \widetilde{V}_k
$ be the inverse branch. In local coordinates on $\widetilde{U}_j$ and $\widetilde{V}_k$, define
$J_k(x)=J_f(\varphi_K)$, $x\in \widetilde{U}_j$. 
The boundedness of the second derivative implies that there are constants
$0<c_k\leq C_k<\infty$ such that
$c_k\leq J_k(x)\leq C_k$. Let \(u\in C^\infty(M)\) be supported in \(F_k^{-1}(\widetilde{U}_j)\), and put \(u_k=u|_{\widetilde{V}_k}\).  Define
$v_k=u_k\circ \varphi_k$. 
Then, by the definition, for \(x\in \widetilde{U}_j\),
\[
   \mathcal{L} u(x)=\sum_{k\in \Theta_j} \frac{u_k(\varphi_k(x))}{J_F(\varphi_k(x))}
        =\sum_{k\in \Theta_j} \frac{v_k(x)}{J_k(x)}.
\]
Thus, for the branch operator \(\mathcal{L}_k\) defined by $\mathcal{L}_k u_k=J_k^{-1}v_k
$, we get
\begin{equation*}\label{LocalL}
   \mathcal{L}u=\sum_{k\in \Theta_j} \mathcal{L}_k u_k.
\end{equation*}
For \(\mu>0\) define two scalar symbols
\[
   \alpha_k(x,\xi)
   =\frac{W^\mu(x,\xi)}{J_k(x)^{1/2}\, W^\mu(DF_k^*(x,\xi))}
,\,\,   \beta_k(x,\xi)
   =\frac{W^\mu(DF_k^*(x,\xi))}{J_k(x)^{1/2}}.
\]
and put $A_k=\op(\alpha_k)$ and $B_k=\op(\beta_k)$. It is straightforward to verify that $A_k \in \Psi^0$ and $B_k \in \Psi^\mu$. 
By the boundedness theorem, $A_k : L^2 \to L^2$ and $B_k : H^\mu \to L^2$
are bounded.\footnote{More generally,
$A_k : H^s \to H^s$ and $B_k : H^s \to H^{s-\mu}$ are bounded for every \(s \in \mathbb{R}\).}
\begin{proposition}\label{BasicPro}
\begin{equation*}
\label{eq:first-local-reduction}
   \big\|\op(W^\mu)\mathcal{L}u\big\|_{L^2}
   \leq \big\|\sum_{k\in \Theta_j} A_k B_k\big\|_{L^2}
      +C\sum_{k\in \Theta
      _j} \|v_k\|_{H^{\mu-1}}.
\end{equation*}
\end{proposition}
Before proving Proposition \ref{BasicPro}, we prove a strightforward lemma. 
\begin{lemma}[Operator Lagrange inequality]\label{MainLemma}
There exists \(C>0\) such that, for every $v_k$, $k\in \Theta_j$, with \(b_k=B_kv_k\),
\begin{equation}\label{Lagrange}
   \sum_{k,l\in \Theta_j}\|A_kb_l\|_{L^2}^2
   -\big\|\sum_{k\in \Theta_j} A_kb_k\big\|_{L^2}^2
   \geq
   \frac12\sum_{k\ne l\in\Theta_j}\|A_kb_l-A_lb_k\|_{L^2}^2
   -C\sum_{k\in \Theta_j}\|v_k\|_{H^{\mu-1/2}}^2.
\end{equation}
\end{lemma}
\begin{proof}
Recall the  scalar Lagrange identity that  says, for complex numbers \(x_1,\dots,x_N\) and \(y_1,\dots,y_N\),
\[
   \sum_{k,l=1}^N |x_k|^2|y_l|^2
   -\big|\sum_{k=1}^N \overline{x_k}y_k\big|^2
   =\frac12\sum_{k\ne l}|\overline{x_k}y_l-\overline{x_l}y_k|^2.
\]
To prove (\ref{Lagrange}), we expand both sides. The left hand side of (\ref{Lagrange}) equals to 
\[
\sum_{k,l}\langle A_kb_l,A_kb_l\rangle
     -\big\langle \sum_kA_kb_k,\sum_lA_lb_l\big\rangle 
   =\sum_{k,l}\langle b_l,A_k^*A_kb_l\rangle
     -\sum_{k,l}\langle A_kb_k,A_lb_l\rangle.
\]
The first item in the right hand equals to 
\[
\frac12\sum_{k\ne l}\|A_lb_k-A_kb_l\|_{L^2}^2 
   =\frac12\sum_{k\ne l}
       \big(
          \|A_lb_k\|_{L^2}^2+
          \|A_kb_l\|_{L^2}^2
          -2\re\langle A_lb_k,A_kb_l\rangle
       \big).
\]
The first two sums are equal after swapping \(k\) and \(l\), so the right hand side of (\ref{Lagrange}) equals to 
\[
  \sum_{k\ne l}\|A_kb_l\|_{L^2}^2
     -\sum_{k\ne l}\re\langle A_lb_k,A_kb_l\rangle.
\]
Now, each \(A_k\) has real scalar symbol \(\alpha_k\in S^0\), thus by using \ref{BPDO}(4)   $A_k^*-A_k\in\Psi^{-1}$ and $[A_k,A_l]\in\Psi^{-1}$.  All discrepancies between the exact scalar identity and the operator identity are finite sums of terms of the form
$
   \langle Rb_k,b_l\rangle
$
with $R:H^{-1/2}\to H^{1/2}$ 
and hence, by duality,
we obtain
\[
   |\langle Rb_k,b_l\rangle|
   \leq C\,\bigl(\|b_k\|_{H^{-1/2}}^2+\|b_l\|_{H^{-1/2}}^2\bigr).
\]
Now, since \(B_k\in\Psi^\mu\),
\[
   \|b_k\|_{H^{-1/2}}
   =\|B_kv_k\|_{H^{-1/2}}
   \leq C\|v_k\|_{H^{\mu-1/2}}.
\]
Summing over the finitely many pairs \((k,l)\) gives the proof.
\end{proof}

\begin{proof}[Proof of Proposition \ref{BasicPro}]
By Lemma \ref{BPDO}(1), we know that $A_kB_k=\op(\alpha_k\beta_k)+R_k$, for some $R_k\in\Psi^{\mu-1}$ and thus $A_kB_k=\op\big(\frac{W^\mu}{J_k}\big)+R_k$. 
On the other hand,
$\op(W^\mu)\op(J_k^{-1})=\op\big(\frac{W^\mu}{J_k}\big)+R_k'$, for some 
$R_l'\in\Psi^{\mu-1}$. Subtracting the last two identities gives
\[
   \op(W^\mu)\op(J_k^{-1})-A_kB_k\in\Psi^{\mu-1}.
\]
Hence there exists \(C>0\), independent of the function to which the operator is applied, such that
\[
   \big\|\bigl(\op(W^\mu)\op(J_k^{-1})-A_kB_k\bigr)w\big\|_{L^2}
   \leq C\|w\|_{H^{\mu-1}}.
\]
Now, put $\mathcal{L}_k u_k=J_k^{-1}v_k=\op(J_k^{-1})v_k.$
We have
\[\op(W^\mu)\mathcal{L}_k u_k
   =\op(W^\mu)\op(J_k^{-1})v_l=A_kB_kv_l+E_kv_k,
\]
where \(E_k\in\Psi^{\mu-1}\).  Using the boundedness theorem, $\|E_kv_k\|_{L^2}\leq C\|v_k\|_{H^{\mu-1}}.$
Summing over the branches gives
\[\|\op(W^\mu)\mathcal{L}u\|_{L^2}\leq 
   \big\|\sum_{k\in \Theta_j} \op(W_\mu)\mathcal{L}_ku_k\big\|_{L^2}   \leq \big\|\sum_{k \in \Theta_j} A_kB_lv_k\big\|_{L^2}
      +C\sum_{k \in \Theta_j} \|v_k\|_{H^{\mu-1}}.
\]
Putting \(b_k=B_kv_k\), this becomes
\begin{equation*}
\label{eq:first-local-reduction}
   \|\op(W^{\mu})\mathcal{L}u\|_{L^2}
   \leq \big\|\sum_{k \in \Theta_j} A_kb_k\big\|_{L^2}
      +C\sum_{k \in \Theta_j} \|v_k\|_{H^{\mu-1}}.
\end{equation*}
\end{proof}
Note that by Lemma \ref{MainLemma}, the following simpler estimate holds
\begin{equation*}
\label{eq:l2-S-bound-by-double-sum}
   \big\|\sum_{k \in \Theta_j}A_kb_k\big\|_{L^2}^2
   \leq
   \sum_{k,l \in \Theta_j}\|A_kb_l\|_{L^2}^2
   +C\sum_{k \in \Theta_l}\|v_k\|_{H^{\mu-1/2}}^2.
\end{equation*}
For fixed \(l\),
\begin{equation}\label{EqMid}
   \sum_{k\in \Theta_j}\|A_kb_l\|_{L^2}^2
   =\sum_{k\in \Theta_l}\langle A_kb_l,A_kb_l\rangle
  =\sum_{k\in \Theta_j}\langle b_l,A_k^*A_kb_l\rangle.
\end{equation}

\medskip
\hspace{-0.4cm}\textbf{Product and adjoint.} Since \(A_k=\op(\alpha_k)\), \(\alpha_k\) is real, and by Lemma \ref{BPDO}(3) we have \(A_k^*A_k=\op(\alpha_k^2)+\Psi^{-1}\), putting 
$a(x,\xi)=\sum_{k \in \Theta_j}\alpha_k(x,\xi)^2$, we get
\[
   \sum_{k \in \Theta_j} A_k^*A_k=\op(a)+R,
   \]
for $R\in\Psi^{-1}$. Consequently, the same \(H^{-1/2}\)-argument as above gives
\[
   |\langle b_l,Rb_l\rangle|
   \leq C\|b_l\|_{H^{-1/2}}^2
   \leq C\|v_l\|_{H^{\mu-1/2}}^2.
\]
Therefore, from (\ref{EqMid}), we get
\begin{equation*}
\label{eq:double-sum-symbol-a}
   \sum_{k,l\in \Theta_j}\|A_kb_l\|_{L^2}^2
   \leq
   \sum_{l\in \Theta_j} \re\langle b_l,\op(a)b_l\rangle
   +C\sum_{l \in \Theta_j}\|v_l\|_{H^{\mu-1/2}}^2.
\end{equation*}
Note that since $\alpha_k(x,\xi)^2=\frac{W^{2\mu}(x,\xi)}{J_k(x)W^{2\mu}(DF_k^*(x,\xi))}$, we have 
\begin{equation}
\label{eq:symbol-a}
   a(x,\xi)=\sum_{k \in \Theta_j}
   \frac{W^{2\mu}(x,\xi)}{J_k(x)W^{2\mu}(DF_k^*(x,\xi))}.
\end{equation}

\medskip
\hspace{-0.4cm}\textbf{Comparison with $B^{2\mu}(F)$.} Now the time has come to incorporate the constants associated with that concept virtually expanding into the calculations. Actually, we show that the principal symbol \(a(x,\xi)\) is bounded from above by \(B^{2\mu}(F)+\eps\) outside a sufficiently large compact set in the cotangent fibers. For fixed \((x,\xi)\), let $ p_k=\varphi_k(x)$ and $\eta_k=(F_k)^*\xi.$
 Then, $DF_k^*(x,\xi)=(p_k,\eta_k).$
\[
   b^{2\mu}((x,\xi);F)
   =\sum_{k\in \Theta_j}
     \frac{(|\xi|_x/|\eta_k|_{p_k})^{2\mu}}{J_F(p_k)} \\
   =\sum_{k\in\Theta_l}
     \frac{|\xi|_x^{2\mu}}{J_k(x)|\eta_k|_{p_k}^{2\mu}}.
\]
On the other hand, by \eqref{eq:symbol-a},
\[
a(x,\xi)=\sum_{k\in \Theta_j}
\frac{\ang{|\xi|_x}^{2\mu}}{J_k(x)\ang{|\eta_k|_{p_k}}^{2\mu}}.
\]
Using the Uniform ellipticity, for every \(\eps>0\), there exists \(R>0\) such that for \(|\xi|_x\geq R\),
\[
   a(x,\xi)
   \leq b^{2\mu}((x,\xi);F)+\eps
   \leq B^{2\mu}(F)+\eps.
\]
Assembling, it gets that for \(|\xi|_x\geq R\), 
\begin{equation*}
\label{eq:a-high-frequency-bound}
   a(x,\xi)
   \leq B^{2\mu}(F)+\eps.   
\end{equation*}
Choose a real number \(\kappa\) satisfying
$\kappa\leq -\frac12$ and $   \kappa<\min\{0,\mu-1/2\}$. Since \(a(x,\xi)\leq B^{2\mu}(F)+\eps\), for \(|\xi|\geq R\), the function $ B^{2\mu}(F)+\eps-a(x,\xi)$
is non-negative outside a compact set in \(T^*U\).  The function \(W^{2\kappa}(x,\xi)\), defined by (\ref{DefW}), is also positive everywhere.  Hence, by choosing \(C_0>0\) sufficiently large, the symbol
\[
   c(x,\xi)
   =B^{2\mu}(F)+\eps+C_0W^{2\kappa}(x,\xi)-a(x,\xi)
\]
is non-negative for every \((x,\xi)\).
Now, since \(a\in S_{1,0}^0\) and \(W^{2\kappa}\in S_{1,0}^{2\kappa}\subset S_{1,0}^0\), so we can apply the Sharp  G\aa rding inequality to the symbol  \(c\) giving that
\[
   Re\langle \op(c)b_k,b_k\rangle
   \geq -C\|b_k\|_{H^{-1/2}}^2.
\]
Now, we use the following three estimates. 
\begin{enumerate}
   
\item The operator \(\operatorname{Op}(W^{2\kappa})\) has order \(2\kappa\).  Its quadratic form is bounded by
$$   |\langle \op(W^{2\kappa})b_k,b_k\rangle|\leq C\|b_k\|_{H^\kappa}^2.$$
\item  Since \(B_k\in\Psi^\mu\), we have $\|b_k\|_{H^\kappa}=\|B_kv_k\|_{H^\kappa}\leq C\|v_k\|_{H^{\mu+\kappa}}.$
\item Since \(\kappa\leq-1/2\) and $\mu+\kappa\leq\mu-1/2$, one gets
\[\|v_k\|_{H^{\mu+\kappa}}
   \leq C\|v_k\|_{H^{\mu-1/2}},
\,\,
\|b_k\|_{H^{-1/2}}
\leq C\|v_k\|_{H^{\mu-1/2}}.
\]

\end{enumerate}
Combining the preceding computation, we obtain
\[
 \re\langle \op(a)b_k,b_k\rangle
\leq (B^{2\mu}(F)+\eps)\|b_k\|_{L^2}^2
+C_0\,\re\langle \op(W^{2\kappa})b_k,b_k\rangle 
 +C\|b_k\|_{H^{-1/2}}^2.
\]
Applying the preceding estimates (1), (2), and (3) to the above inequality, we have
\begin{equation*}
\label{eq:garding-result}
   \re\langle \op(a)b_k,b_k\rangle
   \leq
   (B^{2\mu}(F)+\eps)\|b_k\|_{L^2}^2
   +C\|v_k\|_{H^{\mu-1/2}}^2.
\end{equation*}
Combining the obtained estimate, we obtain the following inequality.
\begin{equation}
\label{eq:S-bound-by-b}
\Big\|\sum_{k\in \Theta_j} A_kb_k\Big\|_{L^2}^2
\leq (B^{2\mu}(F)+\eps)\sum_{k\in \Theta_j}\|b_k\|_{L^2}^2   +C\sum_{k\in \Theta_j}\|v_k\|_{H^{\mu-1/2}}^2.
\end{equation}

\medskip
\hspace{-0.4cm}\textbf{Coordinate changes.} In the next step, we compare the  \(\|b_k\|_{L^2}\) and \(\|\op(W^\mu)u_k\|_{L^2}\). To do this, we define the unitary coordinate-change associated to \(F_k:\widetilde{V_k}\to \widetilde{U_j}\) by
\[
   \calU_k w(x)=\frac{w(\varphi_k(x))}{J_k(x)^{1/2}}.
\]
Indeed,
\[
   \|\calU_k w\|_{L^2(\widetilde{U_j})}^2
   =\int_{\widetilde{U_j}} \frac{|w(\varphi_k(x))|^2}{J_k(x)}\,dm(x)    =\int_{\widetilde{V_k}}|w(p)|^2\,dm(p).
\]
By using the coordinate-change formula for pseudodifferential operators(see \cite[Chapter II, \S5]{Tay81}), we have
\begin{equation}
\label{eq:coordinate-change-pdo}
\op(W^\mu\circ DF_k^*)\calU_k
=\calU_k\op(W^\mu)+R_k^{\mathrm{'''}},
\end{equation}
for some $R_k^{'''}\in\Psi^{\mu-1}$. Now, observe that $\beta_k=(W^\mu\circ DF_k^*)J_k^{-1/2}$. Thus the product formula gives
\[B_kv_k
=\op(W^\mu\circ DF_k^*)\bigl(J_k^{-1/2}v_k\bigr)+R_k''v_k,
\]
for some  $R_k''\in\Psi^{\mu-1}$.
Hence, using \eqref{eq:coordinate-change-pdo} and the fact that $J_k^{-1/2}v_j=\calU_k u_k.$, we get
\begin{align*}
   b_k=B_kv_k
   &=\op(W^\mu\circ DF_k^*)\calU_k u_k+R_k''v_k\\&=\calU_k\op(W^\mu)u_k+R_k^{'''}u_k+R_k''v_k.
\end{align*}
Taking \(L^2\)-norms and using that \(\calU_k\) is unitary,
\begin{equation}\label{Unitary}
\|b_k\|_{L^2}
\leq \|\op(W_\mu)u_k\|_{L^2}
 +C\|u_k\|_{H^{\mu-1}}
+C\|v_k\|_{H^{\mu-1}}.
\end{equation}
But the composition with a fixed smooth diffeomorphism preserves Sobolev spaces of every order, so
$   \|v_k\|_{H^{\mu-1}}\leq C\|u_k\|_{H^{\mu-1}}.
$
Therefore (\ref{Unitary}) gives
\begin{equation}
\label{eq:bj-by-uj}
\|b_k\|_{L^2}
\leq \|\op(W^\mu)u_k\|_{L^2}+C\|u_k\|_{H^{\mu-1}}.
\end{equation}
Squaring \eqref{eq:bj-by-uj} gives, for any \(\delta>0\),
\[
\|b_k\|_{L^2}^2
\leq (1+\delta)\|\op(W^\mu)u_k\|_{L^2}^2
 +C_\delta\|u_k\|_{H^{\mu-1}}^2.
\]
By choosing \(\delta>0\) small enough such that
$(B^{2\mu}(F)+\eps/2)(1+\delta)
   \leq B^{2\mu}(F)+\eps.
$, \eqref{eq:S-bound-by-b} implies
\begin{equation}
\label{eq:S-bound-by-u}
   \Big\|\sum_{k\in \Theta_j} A_kb_k\Big\|_{L^2}^2
   \leq
   (B^{2\mu}(F)+\eps)
   \sum_{k\in \Theta_j}\|\op(W^\mu)u_k\|_{L^2}^2
   +C\sum_{k\in \Theta_j}\|u_k\|_{H^{\mu-1/2}}^2.
\end{equation}

\medskip
\begin{proof} [Proof of Local Lasota-York inequality]Let \(\mu^\prime\) satisfies $\max\{0,\mu-1/2\}\leq \mu'<\mu.$ Then $$\|u_k\|_{H^{\mu-1/2}}
   \leq C\|u_k\|_{H^{\mu'}}.$$
Furthermore, since \(\mu-1\leq\mu-1/2\leq\mu'\), we also have $ \|u_k\|_{H^{\mu-1}}
   \leq C\|u_k\|_{H^{\mu'}}.$
Composition with \(\varphi_k\) gives the same bounds for \(v_k\). That is
\begin{equation}\label{LocalLY}
   \|v_k\|_{H^{\mu-1}}
   +\|v_k\|_{H^{\mu-1/2}}
   \leq C\|u_k\|_{H^{\mu'}}.
\end{equation}
Combining \eqref{LocalLY} with \eqref{eq:first-local-reduction} and \eqref{eq:S-bound-by-u}, we get
\begin{equation}
\label{eq:local-LY}
\|\op(W^\mu)\mathcal{L}u\|_{L^2}
\leq
\sqrt{B^{2\mu}(F)+\epsilon}
\Big(\sum_{k\in \Theta_j}\|\op(W^\mu)u_k\|_{L^2}^2\Big)^{1/2}
+C\sum_{k\in \Theta_j}\|u_k\|_{H^{\mu'}}.
\end{equation}
Now, \eqref{Sineq}, \eqref{LocalL} and \eqref{eq:local-LY} give the desired inequality
$$\|\mathcal{L}u\|\leq 
\sqrt{B^{2\mu}(F)+\epsilon}\|u\|_{H^\mu}+
C\|u\|_{H^{\mu'}}
.$$
\end{proof}
\section*{Acknowledgements}
The authors is deeply grateful to Masato Tsujii for his invaluable guidance and 
insightful suggestions and detailed feedback during our discussions.
The first author is partially supported by IPM grant NO.1404340211.

\bibliography{ref.bib}
\bibliographystyle{amsplain}

\end{document}